\documentclass[10pt]{article}
\usepackage{latexsym,amsmath,amssymb,graphics,amscd, empheq}
\usepackage{enumitem}
\usepackage{xcolor}
\usepackage{stackengine}
\usepackage[breaklinks,colorlinks,
  citecolor=blue,
  urlcolor=blue]{hyperref}
\DeclareFontFamily{OT1}{pzc}{}
\DeclareFontShape{OT1}{pzc}{m}{it}{<-> s * [1.10] pzcmi7t}{}
\DeclareMathAlphabet{\mathpzc}{OT1}{pzc}{m}{it}

\addtocounter{footnote}{1}
\makeatletter
\@addtoreset{table}{bsection}
\def\thetable{\thesection.\@arabic\c@table}
\def\fps@table{h, t}
\@addtoreset{equation}{section}

\makeatother

\newtheorem{theorem}{Theorem}[section]
\newtheorem{definition}[theorem]{Definition}
\newtheorem{lemma}[theorem]{Lemma}
\newtheorem{remark}[theorem]{Remark}
\newtheorem{proposition}[theorem]{Proposition}

\newtheorem{corollary}[theorem]{Corollary}
\newtheorem{example}[theorem]{Example}

\newcommand{\bfi}{\bfseries\itshape}
\newcommand{\vertiii}[1]{{\left\vert\kern-0.25ex\left\vert\kern-0.25ex\left\vert #1 
    \right\vert\kern-0.25ex\right\vert\kern-0.25ex\right\vert}}
\newsavebox{\savepar}

\makeatletter

\usepackage{scalerel,stackengine}
\stackMath
\newcommand\reallywidehat[1]{%
\savestack{\tmpbox}{\stretchto{%
  \scaleto{%
    \scalerel*[\widthof{\ensuremath{#1}}]{\kern-.6pt\bigwedge\kern-.6pt}%
    {\rule[-\textheight/2]{1ex}{\textheight}}%WIDTH-LIMITED BIG WEDGE
  }{\textheight}% 
}{0.5ex}}%
\stackon[1pt]{#1}{\tmpbox}%
}
\begin{document}

\title{\textbf{Dimension reduction in  recurrent networks by canonicalization}}
\author{Lyudmila Grigoryeva$^{1}$, and Juan-Pablo Ortega$^{2}$}
\date{}
\maketitle

\begin{abstract}
Many recurrent neural network machine learning paradigms can be formulated using state-space representations. The classical notion of canonical state-space realization is adapted in this paper to accommodate semi-infinite inputs so that it can be used as a dimension reduction tool in the recurrent networks setup. The so-called input forgetting property is identified as the key hypothesis that guarantees the existence and uniqueness (up to system isomorphisms) of canonical realizations for causal and time-invariant input/output systems with semi-infinite inputs. Additionally, the notion of optimal reduction coming from the theory of symmetric Hamiltonian systems is implemented in our setup  to construct canonical realizations out of input forgetting but not necessarily canonical ones. These two procedures are studied in detail in the framework of linear fading memory input/output systems. Finally, the notion of  implicit reduction using reproducing kernel Hilbert spaces (RKHS) is introduced which allows, for systems with linear readouts, to achieve dimension reduction without the need to actually compute the reduced spaces introduced in the first part of the paper.
\end{abstract}

\bigskip

\textbf{Key Words:} recurrent neural network, reservoir computing, dimension reduction, state-space system, canonicalization, echo state network, ESN, linear recurrent network, machine learning, echo state property.

\makeatletter
\addtocounter{footnote}{1} \footnotetext{%
Department of Statistics, University of Warwick, Coventry CV4 7AL , UK.}
\makeatother
\makeatletter
\addtocounter{footnote}{1} \footnotetext{%
Division of Mathematical Sciences, 
Nanyang Technological University,
21 Nanyang Link,
Singapore 637371.
{\texttt{Juan-Pablo.Ortega@ntu.edu.sg}}}
\makeatother

\medskip

\medskip

\medskip

\section{Introduction}

State-space models are of widespread use in the construction of input/output systems in many application contexts. The Markovian nature of the state equation makes them particularly convenient in the construction of efficient simulation algorithms without preventing the possibility of encoding long-memory type behaviors. These models were first introduced in the context of systems and control theory \cite{kalman1959general, kalman1959unified, kalman:original, Kalman1962, baum1966statistical, Kalman2010} and met spectacular success in all sorts of industrial, military, and scientific applications in relation to filtering, smoothing, and forecasting (see \cite{kalman1960new, kalman1961new, hutchinson1984kalman, Koopman:kalman, Sarkka2013} and references therein for just a few examples). 

More recently, these systems have reemerged in the context of the machine learning of dynamic processes as powerful recurrent network paradigms. The question of interest in this framework is the learning or the estimation of the parameters of a state-space system out of finite-length realizations of the input and output processes. This learning problem is, to some extent, just a reformulation of the non-linear identification problem that has been thoroughly studied in systems and control theory \cite{kalman:mit, Sontag1979, DangVanMien1984, Narendra1990, Matthews1994, lindquist:picci} as well as in the theory of empirical processes \cite{Marie:Duflo}. 

Despite these similarities, there are new problems that need to be seriously addressed when using state-space systems in the machine learning context. For instance, much of the systems theory literature is dedicated to the characterization of the controllability question for invertible systems and  formulated using a prescribed initial or final condition (see \cite{fliess1981group, normand1983theorie, jakubczyk1990controllability} for an in-depth study of the discrete-time case). However, in most machine learning situations, it is more appropriate to work using {\bfi  semi-infinite temporal traces towards the past} in which the dependence on initial conditions disappears. This feature arises in the presence of time-invariant input/output systems and stationary stochastic processes, and it is a crucial element in the formulation of the fading memory property that pervades many modeling situations. Additionally, most systems that are considered in applications are subsystems of a Markovian system that, generically, exhibit a functional dependence on the infinite past. 

Another distinctive feature of state-space models in the learning framework is the use of  randomization. Early in the application of these models as recurrent networks, important difficulties were identified at the time of their training using  classical gradient descent (backpropagation-type) methods having to do with bifurcation phenomena~\cite{Doya92} in these intrinsic dynamical models. 
Recent progress in the regularization and training of recurrent structures (see, for instance \cite{Graves2013, pascanu:rnn, zaremba}, and references therein) solves to some extent some of these non-convergence problems. A different approach to circumvent this question, specially in data-intensive applications, is to use randomly generated state equations and to only train the time-independent observation equation that is selected out of a functionally simple (preferably linear) family. This revolutionary idea has its origin in static frameworks like, for instance, in the seminal works on random feature models \cite{Rahimi2007} and Extreme Learning Machines \cite{Huang2006}. This philosophy was extended to the dynamical context that we are interested in this paper under the names of {\bfi  reservoir computing (RC)} \cite{jaeger2001, Jaeger:2002, Jaeger04} and {\bfi  liquid state machines} \cite{maass1, maass2} and has proved to be very successful in a great variety of empirical classification and forecasting applications (see, for instance, \cite{Jaeger04, Wyffels2008, lukosevicius, Wyffels2010, Buteneers2013, GHLO2014, Ott2018, Pathak:PRL, pathak2018}).

These empirical discoveries have motivated an intense activity in the theoretical front to understand, quantify, and optimize the information processing abilities of state-space systems. An important body of work has to do with the  assessment of the memory and forecasting abilities of these constructions in terms of their architectures and dependence properties of the input signals~\cite{Jaeger:2002, White2004, Ganguli2008, Hermans2010, dambre2012, esn2014, GHLO2014, GHLO2014_capacity, linearESN, RC4pv, farkas:bosak:2016, Goudarzi2016, RC3, Xue2017, Charles2017, marzen:capacity, Verzelli2019a, RC15}. Additionally, memory capacities have been extensively compared with other related concepts like Fisher information-based criteria \cite{Tino2013, Livi2016, tino:symmetric}.

In a more learning theoretical note, much progress has been done in the last years in the understanding of the universal approximation and the generalization properties of this approach. By now, we can find in the literature many families of state-space systems that have been proved to be universal approximants in different contexts. For example, when inputs are deterministic and uniformly bounded, universality has been proved for linear systems with polynomial observation equations \cite{Boyd1985, RC6},  state-affine systems (SAS) \cite{RC6},  the echo state networks (ESNs) \cite{RC7, RC20} introduced in \cite{Matthews:thesis, Matthews1994, Jaeger04}, the so-called signature state-affine systems (SigSAS) \cite{RC13} that encode in state-space form the truncation of Volterra series expansions, or the temporal convolutional networks \cite{Hanson2019}. These results have been extended to a stochastic setup in \cite{RC8} and also exist in the context of the approximation of dynamical systems with a compact phase space \cite{hart:ESNs, allen:tikhonov, RC18, RC21}. By now, risk \cite{RC10} and approximation \cite{RC12} bounds exist for some of these systems similar to those that can be formulated, for instance, in the context of shallow neural networks or other static machine learning paradigms.

In this paper we focus on another machine learning aspect of major importance in the practical use of reservoir computing and state-space systems, namely, {\bfi  dimension reduction}. Given a machine learning paradigm, the dimension reduction problem consists of generically  finding a system with reduced complexity that exhibits equivalent or almost equivalent approximation properties.  For example, in the feedforward neural networks context, there exist standard pruning techniques \cite{Haykin2009} that determine which neurons can be eliminated in a given network configuration when they are not relevant for a given approximation task. Other widespread strategies consist in using principal components analysis or random projections in the spirit of \cite{JLlemma} (see \cite{RC13} for a first step in the use of these techniques in reservoir computing). 

In the framework of mechanical and controlled systems, dimension reduction is a classical and well-studied subject that goes back to Jacobi's elimination of the node in multi-body celestial mechanics in the nineteenth century. In that setup, dimension reduction is, most of the time, associated with the use of the conserved quantities associated to the symmetries of a given system and that are encoded in the level sets of a momentum map \cite{kostant66, souriau1966, souriau, smale}. Dimension reduction is generically obtained by restricting the dynamics to invariant manifolds and by projecting it onto the orbit space with respect to the residual symmetry that leaves those invariant. In the context of autonomous systems, this procedure is referred to as {\bfi  Marsden-Weinstein reduction} \cite{mwr}; see \cite{Ortega2004, Marsden2007} for self-contained presentations of this beautiful theory. Part of these mostly differential geometric techniques for dimension reduction has been extended to controlled systems. See, for instance, \cite{van1981symmetries, nijmeijer1982controlled, grizzle1985structure, van1987symmetries, blankenstein2004singular, gay2011clebsch, ohsawa2013symmetry, bloch2015nonholonomicbook} and references therein.

Many reservoir computing applications like, for instance, those in \cite{Jaeger04, Ott2018, Pathak:PRL, pathak2018} require the use of systems with state-space dimensions in the thousands that, generically, present no symmetries that could be used for reduction. This motivates the investigation of  another natural dimension reduction related notion, this time only applicable to state-space systems, namely that of {\bfi  canonicalization}. The idea behind it is based on the observation that since the state-space representation of input/output systems is not unique, one should choose the most ``economical" one in which ``unused" states are dropped from the representation and those that are ``undistinguishable" from a dynamical point of view are identified by the passage to a quotient space. These ``optimal" state-space representations are called {\bfi  canonical realizations}, and in the context of forward looking systems it can be proved that they exist and are unique up to system isomorphisms. This result is usually called the {\bfi  Canonical Realization Theorem} (see for instance \cite[Chapter 2]{Matthews:thesis}).

The main goal of this paper is {\it extending these canonicalization results to the context of time-invariant and causal input/output systems with semi-infinite inputs} and, moreover, to obtain a Canonical Realization Theorem in this framework out of a reduction approach similar to the one introduced in \cite{optimal:mm, reduction:optimal:cras}. More explicitly, the paper contains two main canonicalization results:
\begin{itemize}
\item A {\bfi  Canonical Realization Theorem} (Theorem \ref{Canonical realization of input/output systems}) for input/output systems. This result shows that any {\bfi  causal and time-invariant} filter that has the so-called {\bfi  input forgetting property}  admits a canonical state-space realization that is unique up to system isomorphisms. The input forgetting property (also referred to in the literature as the {\bfi   unique steady-state property}) is a modeling feature that appears profusely in applications and that can be obtained out of the so-called {\bfi  fading  memory property} (see \cite{Boyd1985, RC9} for a detailed discussion about these concepts). An important merit of Theorem \ref{Canonical realization of input/output systems} is identifying the input forgetting property as the key concept that leads to the availability of canonical realizations in the presence of semi-infinite inputs. Additionally, it constitutes a result of great generality as it provides a constructive procedure for the design of state-space realizations for a vast category of  input/output systems; the price to pay for this generality is the potentially complicated nature of the representing state space or its infinite dimensional character (when such notion is well-defined).
\item A  {\bfi  Canonicalization by Reduction Theorem} (Theorem \ref{Canonicalization by reduction}). This result uses a reduction approach similar to the one introduced in \cite{optimal:mm, reduction:optimal:cras} in the context of symmetric Hamiltonian systems to construct a canonical realization for a state-space system that has the input forgetting property system by using a ``reduced" version of it in a sense that will be introduced in detail later on.
\end{itemize}

These two results are illustrated and applied in detail in Section \ref{Realization and canonicalization of linear filters} in the context of linear fading memory filters. In particular, Theorem \ref{Canonicalization of linear filters} shows that any linear, causal, time-invariant filter with semi-infinite inputs that has the fading memory property (that, as we shall see, implies the input forgetting property) admits a canonical linear state-space realization (possibly infinite dimensional). Additionally, this result also characterizes all the isomorphic canonical realizations of the given filter as a homogeneous manifold constructed using the general linear group of the state space. Finally, the  Canonicalization by Reduction Theorem  \ref{Canonicalization by reduction} in the linear setup yields Theorem \ref{Canonicalization by reduction of linear state-space systems}, which fully characterizes how to construct a canonical linear realization by shrinking  the linear state-space appropriately, for a given linear system that has the input forgetting property but that is not necessarily canonical. 

{The paper concludes with Section \ref{Implicit reduction using RKHS}, where we introduce what we call {\bfi  implicit reduction} using reproducing kernel Hilbert spaces (RKHS). The main goal of that section consists in circumventing the need of computing the reduced spaces introduced in the previous sections, which may be technically difficult, in order to achieve dimension reduction. As we show in those pages, the RKHS formulation of the estimation problem for state-space systems with linear readouts achieves exactly that as a consequence of the well-known Representer Theorem \cite[page 117]{Mohri:learning:2012}. Section \ref{Conclusions} concludes the paper.}

\section{Canonical systems with semi-infinite inputs}
\label{Recurrent neural networks with stationary inputs}

We briefly introduce a few definitions that make explicit the setup where we shall be working. The objects of interest in this paper are input/output systems determined by state-space systems. The symbols ${\cal Z}  $  and ${\cal Y} $ will denote the {\bfi  input} and the {\bfi  output spaces}, respectively, and ${\cal X} $ will be the {\bfi  state space} of the system that will create the link between them. These three spaces are typically subsets of a Euclidean space or, more generally, finite or infinite dimensional manifolds; for the time being we shall assume no particular structure on them. A {\bfi  discrete-time state-space system}  is determined by the following two equations that put in relation sequences ${\bf z} \in {\cal Z}^{\mathbb{Z}}, {\bf y} \in {\cal Y}^{\mathbb{Z}}, {\bf x} \in {\cal X}^{\mathbb{Z}} $ in the three spaces that we just introduced:
\begin{empheq}[left={\empheqlbrace}]{align}
\mathbf{x} _t &=F(\mathbf{x} _{t-1}, {\bf z}_t),\label{rc state eq}\\%\sigma \left(A\mathbf{x}_{t-1}+ C{ z} _t+ \boldsymbol{\zeta}\right),\label{esn reservoir equation theorem prep}\\
{\bf  y} _t &= h( \mathbf{x} _t) , \label{rc readout eq}
\end{empheq}
for any $t \in \Bbb Z $. The map $F:  {\cal X}  \times {\cal Z} \longrightarrow {\cal X}$ is called the {\bfi  state map} and  $h:{\cal X} \longrightarrow {\cal Y}$ the {\bfi  readout } or {\bfi  observation } map. We shall sometimes denote a system by using the  triple $({\cal X}, F, h)$. The term {\bfi  recurrent neural network (RNN)}  is used sometimes in the literature to refer to state-space systems where the state map $F$ in \eqref{rc state eq} is neural network-like, that is, it is a concatenation of compositions of a nonlinear activation function with an affine function of the states and the input. A particular case of RNNs are the {\bfi  echo state networks} introduced in \cite{Matthews:thesis, Jaeger04} where one neural layer of this type is used (with random connectivity between neurons in \cite{Jaeger04}).

We  focus on state-space systems of the type \eqref{rc state eq}-\eqref{rc readout eq} that determine an {\bfi  input/output} system. This happens in the presence of the so-called {\bfi  echo state property (ESP)}, that is, when for any  ${\bf z} \in {\cal Z}^{\mathbb{Z}}$ there exists a unique $\mathbf{y} \in {\cal Y}^{\mathbb{Z}}$ such that \eqref{rc state eq}-\eqref{rc readout eq} hold. In that case, we  talk about the {\bfi  state-space filter} $U ^F_h: {\cal Z}^{\mathbb{Z}}\longrightarrow  {\cal Y}^{\mathbb{Z}} $ associated to the state-space system $({\cal X}, F, h)$  defined by: 
\begin{equation*}
U  ^F _h({\bf z}):= \mathbf{y},
\end{equation*}
where ${\bf z} \in {\cal Z}^{\mathbb{Z}} $  and $\mathbf{y} \in {\cal Y}^{\mathbb{Z}} $ are linked by \eqref{rc state eq}-\eqref{rc readout eq} via the ESP. If the ESP holds at the level of the state equation \eqref{rc state eq}, we can define a {\bfi  state filter}  $U ^F: {\cal Z}^{\mathbb{Z}}\longrightarrow  {\cal X}^{\mathbb{Z}} $ and, in that case, we have that
\begin{equation*}
U ^F_h:= h \circ U ^F. 
\end{equation*}
It is easy to show that state and state-space filters are automatically causal and time-invariant (see \cite[Proposition 2.1]{RC7}) and hence it suffices to work with their restriction $U ^F_h: {\cal Z}^{\mathbb{Z}_-}\longrightarrow  {\cal Y}^{\mathbb{Z}_-} $ to semi-infinite inputs and outputs. Moreover, $U ^F_h $  determines a state-space {\bfi   functional} $H ^F_h: {\cal Z}^{\mathbb{Z}}\longrightarrow  {\cal Y}$ as $H ^F_h({\bf z}):=U ^F_h({\bf z})_0 $, for all ${\bf z} \in {\cal Z}^{\mathbb{Z}_{-}} $ (the same applies to $U ^F $ and $H ^F $ when the ESP holds at the level of the state equation). In the sequel we use the symbol $\mathbb{Z}_{-}  $ to denote the negative integers including zero and $\mathbb{Z}^{-}  $ without zero.

\paragraph{State-space morphisms.} 
As we already mention in the introduction, a given input/output filter may have different state-space realizations. One way to construct them is by using the natural functors  between state-space systems that we define below. Consider the state-space systems determined by the triples $({\cal X} _i, F _i, h _i)$, $i \in \left\{1,2\right\}$, with $F _i: {\cal X} _i\times {\cal Z}\longrightarrow {\cal X}_i$ and $h _i:{\cal X}_i \longrightarrow {\cal Y} $.
\begin{definition}
A map $f: {\cal X} _1 \longrightarrow {\cal X}_2$ is a {\bfi  morphism} between the systems $({\cal X} _1,F_1, h_1)$ and $({\cal X} _2, F_2, h_2)$ whenever it satisfies the following two properties:
\begin{description}
    \item[(i)] {\bfi System equivariance:} $f(F_1({\bf x}_1, {\bf z})) = F_2(f({\bf x}_1), {\bf z})$, for all ${\bf x}_1 \in {\cal X}_1$ and ${\bf z} \in {\cal Z}$.
    \item[(ii)] {\bfi Readout invariance:} $h_1({\bf x}_1) = h_2(f({\bf x}_1))$, for all ${\bf x}_1 \in {\cal X} _1$.
\end{description}
\end{definition}

When the map $f$ has an inverse $f ^{-1} $  and this inverse is also a morphism between the systems determined by  $({\cal X} _2, F_2, h_2)$ and $({\cal X} _1, F_1, h_1)$ and  we say that $f$ is a {\bfi system isomorphism} and that the systems $({\cal X} _1, F_1, h_1)$ and $({\cal X} _2, F_2, h_2)$ are {\bfi isomorphic}. We note that given a system $F_1: {\cal X}_1 \times {\cal Z} \longrightarrow {\cal X}_1, h_1: {\cal X}_1 \longrightarrow {\cal Y}$ and a bijection $f:{\cal X}_1 \longrightarrow {\cal X}_2$, the map $f$ is a system isomorphism with respect to the system $F_2: {\cal X}_2 \times {\cal Z} \longrightarrow {\cal X}_2, h_2: {\cal X}_2 \longrightarrow {\cal Y}$ defined by 
\begin{align}
    F_2({\bf x}_2, {\bf z}) &:= f(F_1(f^{-1}({\bf x}_2), {\bf z})), \quad \text{for all} \quad {\bf x}_2 \in {\cal X}_2, {\bf z} \in {\cal Z}, \label{isomorphic state map}\\
    h_2({\bf x}_2) &:= h_1(f^{-1}({\bf x}_2)), \quad \text{for all} \quad {\bf x}_2 \in {\cal X}_2.\label{isomorphic readout map}
\end{align}
The proof of the following elementary result can be found in \cite{RC15}.

\begin{proposition}
\label{properties of morphisms solutions}
Let  $({\cal X}_i, F _i, h _i)$, $i \in \left\{1,2\right\}$, be two systems with $F _i: {\cal X}_i\times {\cal Z}\longrightarrow {\cal X}_i$ and $h _i:{\cal X}_i \longrightarrow {\cal Y}$.  Let $f: {\cal X}_1 \longrightarrow {\cal X}_2$ be a map. Then:
\begin{description}
\item[(i)] If $f$ is system equivariant and ${\bf x}^1 \in {\cal X}_1^{\mathbb{Z}_-}$ is a solution for the state system associated to $F_1$ and the input ${\bf z} \in {\cal Z}^{\mathbb{Z}_-}$, then so is $(f({\bf x}_t^1))_{t\in \mathbb{Z}_-} \in {\cal X}_2^{\mathbb{Z}_-}$ for the system associated to $F_2$ and the same input.
\item[(ii)] Suppose that the system determined by $({\cal X}_2,F_2, h_2)$ has the  echo state property and assume that the state system determined by $F_1$ has at least one solution for each element ${\bf z} \in {\cal Z}^{\mathbb{Z}_-}$. If $f$ is a morphism between $({\cal X}_1,F_1, h_1)$ and $({\cal X}_2,F_2, h_2)$, then $({\cal X}_1,F_1, h_1)$ has the echo state property and, moreover,
\begin{equation}\label{eq: morphism in proposition}
    U_{h_1}^{F_1} = U_{h_2}^{F_2}.
\end{equation}
\item [(iii)] If $f$ is a system isomorphism, then the implications in the previous two points are reversible, that is, the indices $1$ and $2$ can be exchanged.
\end{description}
\end{proposition}

\paragraph{Reachability and observability.} 
We just showed in the previous paragraph that system morphisms produce different state-space system realizations for a given input/output system. We now introduce dynamical properties that ensure that the reverse implication holds, that is, if we have two different state-space system realizations for a given input/output system we can ensure that there exists a system morphism between them. The following definitions are natural adaptations of the concepts with the same name in the context of forward-in-time systems \cite{sontag:book, Lewis2002, Bullo2005}.

The definition uses the following notation: if ${\bf z} \in {\cal Z} ^{\mathbb{Z}_{-}} $  and $\widetilde{{\bf z}} \in {\cal Z}^T $ for some $T \in \mathbb{N} $,  then the symbol ${\bf z} \widetilde{{\bf z} }\in {\cal Z} ^{\mathbb{Z}_{-}}$ denotes the semi-infinite sequence obtained by concatenation of ${\bf z} $ and $\widetilde{{\bf z}} $. 

\begin{definition}
\label{Reachable and Observable}
Let $({\cal X}, F, h) $ be a state-space system with $F: {\cal X}\times {\cal Z} \longrightarrow {\cal X} $ and $h : {\cal X} \longrightarrow {\cal Y} $. Assume that $({\cal X}, F, h) $ has the echo state property. Then, we say that $({\cal X}, F, h) $ is:
\begin{description}
\item [(i)] {\bf Reachable} (respectively, {\bf strongly reachable}), when for any ${\bf y} \in {\cal Y} $ (respectively, $\mathbf{x} \in {\cal X} $) there exists ${\bf z} \in {\cal Z} ^{\mathbb{Z}_{-}} $ such that $H_h^F ({\bf z})= {\bf y} $ (respectively, $H^F({\bf z})= \mathbf{x} $).
\item [(ii)] {\bf Observable}, when it does not have indistinguishable states. Two {distinct} states $\mathbf{x} _1, \mathbf{x} _2 \in {\cal X}  $ are called {\bfi  indistinguishable} when {there exist} $ {\bf z} _1, {\bf z} _2  \in {\cal Z}^{\mathbb{Z}_{-}}  $ such that $\mathbf{x} _1= H ^F({\bf z} _1)$, $\mathbf{x} _2= H ^F({\bf z} _2)$ {and, additionally,} we have that $H ^F_h({\bf z} _1\widetilde{{\bf z}})=H ^F_h({\bf z} _2\widetilde{{\bf z}})$, for any $\widetilde{{\bf z}} \in \mathcal{Z} ^T$ and any $T \in \mathbb{N} $.
\item [(iii)] {\bf Canonical}, when $({\cal X}, F,h)$ is strongly reachable and observable.
\end{description}
\end{definition}

\noindent Note that if the observation map $h$ is surjective, then strong reachability implies reachability.

\begin{proposition}
\label{proposition existence of morphism}
Let $({\cal X},F,h ) $ and $(\overline{{\cal X}},\overline{F},\overline{h} ) $ be two systems that have the echo state property  and yield the same time-invariant  input/output system, that is, $H ^F_h=H ^{ \overline{F}}_{ \overline{h}} $. If $({\cal X},F,h ) $ is strongly reachable and $(\overline{{\cal X}},\overline{F},\overline{h} ) $ is observable then there exists a unique system morphism $f:\mathcal{X} \longrightarrow \overline{\mathcal{X}} $.
\end{proposition}

Before we proceed with the proof of this proposition, we list in the following lemma three elementary properties of time-invariant state-space filters. In the proof we use the {\bfi  time delay} operators $T _\tau : {\cal Z}^{\mathbb{Z}_{-}}\longrightarrow {\cal Z}^{\mathbb{Z}_{-}} $ that, for any $\tau \in \mathbb{N} $, are defined as 
\begin{equation}
\label{first time delay operator}
T _\tau({\bf z}) _t:= {\bf z}_{t- \tau} , \ t \in \mathbb{Z}_{-}. 
\end{equation}
We recall that filters are called time-invariant when they commute with the time delay operators. {Additionally, we will be using the notion of invertible state map. We recall that the map $F: {\cal X} \times {\cal Z} \longrightarrow {\cal X}$ is {\bfi  invertible} when for any ${\bf z} \in {\cal Z} $, the maps $F (\cdot , {\bf z}): {\cal X} \longrightarrow {\cal X} $ are injective and hence there exists a map $F ^{-1}: {\cal X} \times {\cal Z} \longrightarrow {\cal X} $ such that 
\begin{equation}
\label{definition of finverse}
F^{-1}(F ( \mathbf{x}, {\bf z}), {\bf z})= \mathbf{x}, \quad \mbox{for all} \quad \mathbf{x} \in {\cal X}\text{ and } {\bf z} \in {\cal Z}.
\end{equation}
}

\begin{lemma}
\label{preliminary for canonical}
Let $({\cal X},F,h ) $ be a system that has the echo state property with input and output spaces ${\cal Z}  $ and ${\cal Y} $, respectively,  and $\widetilde{{\bf z} } \in {\cal Z} $. Then 
\begin{equation}
\label{first for lemma}
U ^F({\bf z}\widetilde{{\bf z}})_{-1}=H^F({\bf z} ) \quad \mbox{and} \quad
H^F({\bf z}\widetilde{{\bf z}})=F(H ^F({\bf z}), \widetilde{{\bf z}}),\enspace  {\rm for} \enspace {\bf z}\in  {\cal Z}^{\mathbb{Z}_{-}}  .
\end{equation} 
Additionally, if ${\bf z} _1, {\bf z} _2 \in {\cal Z} ^{\mathbb{Z}_{-}} $ are such that $H ^F({\bf z} _1)=H ^F({\bf z} _2)$ then 
\begin{equation}
\label{second for lemma}
H ^F({\bf z} _1\widetilde{{\bf z} })=H ^F({\bf z} _2\widetilde{{\bf z} }), \quad \mbox{for any} \quad \widetilde{{\bf z}} \in {\cal Z}^T \quad \mbox{and any} \quad T \in \mathbb{N}.
\end{equation}
{The converse holds when $F$ is an invertible state map.}
\end{lemma}

\noindent\textbf{Proof of the Lemma.\ \ } The identities in \eqref{first for lemma} are a consequence of the time-invariance of $U^F $. Indeed,
\begin{equation*}
U ^F({\bf z}\widetilde{{\bf z}})_{-1}= \left(T _1\circ U ^F({\bf z}\widetilde{{\bf z}})\right) _0=\left(U ^F (T_1 ({\bf z}\widetilde{{\bf z}}))\right) _0=U ^F ({\bf z}) _0=H^F({\bf z} ).
\end{equation*}
As to the second equality in \eqref{first for lemma}, by definition and the identity that we just proved:
\begin{equation*}
H^F({\bf z}\widetilde{{\bf z}})=U^F({\bf z}\widetilde{{\bf z}}) _0= F(U^F({\bf z}\widetilde{{\bf z}}) _{-1}, \widetilde{{\bf z}})=F(H ^F({\bf z}), \widetilde{{\bf z}}).
\end{equation*}
Concerning \eqref{second for lemma}, let $\widetilde{{\bf z}}=(\widetilde{{\bf z}} _1, \ldots, \widetilde{{\bf z}} _T)\in {\cal Z} ^T$. Then, by the hypothesis $H ^F({\bf z} _1)=H ^F({\bf z} _2)$ and the identity that we just proved:
\begin{equation*}
H^F( {\bf z}_1\widetilde{{\bf z}} _1)= F(H ^F({\bf z}_1), \widetilde{{\bf z}}_1)= F(H ^F({\bf z}_2), \widetilde{{\bf z}}_1)=H^F( {\bf z}_2\widetilde{{\bf z}} _1).
\end{equation*}
Analogously,
\begin{equation*}
H^F( {\bf z}_1\widetilde{{\bf z}} _1\widetilde{{\bf z}}_2)= F(H ^F({\bf z}_1 \widetilde{{\bf z}}_1), \widetilde{{\bf z}}_2)= F(H ^F({\bf z}_2\widetilde{{\bf z}}_1), \widetilde{{\bf z}}_2)=H^F( {\bf z}_2\widetilde{{\bf z}} _1\widetilde{{\bf z}}_2).
\end{equation*}
Repeating this procedure $T$ times yields \eqref{second for lemma}. {Suppose now that $F$ is invertible and that \eqref{second for lemma} holds. In particular, we have that $H ^F({\bf z} _1\widetilde{{\bf z} })=H ^F({\bf z} _2\widetilde{{\bf z} })$,  for any  $\widetilde{{\bf z}} \in {\cal Z}$ which, by \eqref{first for lemma}, implies that $F(H ^F({\bf z}_1), \widetilde{{\bf z}})= F(H ^F({\bf z}_2), \widetilde{{\bf z}}) $. If we now apply $F ^{-1}$ to both sides of this equality we have by \eqref{definition of finverse} that $H ^F({\bf z}_1)= H ^F({\bf z}_2)$, as required.}\quad $\blacktriangledown$

\medskip

{These facts can be used to prove that any system that has the echo state property at the level of the state equation can be restricted to a smaller state space where it becomes strongly reachable. Additionally, they also imply that invertible state maps and injective readouts determine observable state-space systems.}

\begin{corollary}
\label{ESP implies exists reachable set}
{
Let  $({\cal X},F,h ) $ be a system with input space ${\cal Z}  $ that has the echo state property at the level of the state equation. Then there exists a subset ${\cal X}' \subset {\cal X}$ such that $F$ restricts to a map (denoted with the same symbol) $F: {\cal X}' \times {\cal Z} \longrightarrow {\cal X}' $ and, moreover, $({\cal X}',F,h ) $ is strongly reachable.}

{
Additionally, if the map $F: {\cal X} \times {\cal Z} \longrightarrow {\cal X}$ is invertible, then the system $({\cal X},F,h ) $ is necessarily observable for any  readout map $h: {\cal X} \longrightarrow{\cal Y}$ that is injective when restricted to ${\cal X}':=H^F({\cal Z}^{\mathbb{Z}_{-}})$.
}
\end{corollary}

\medskip

\noindent {\textbf{Proof of the Corollary.\ \ } First, the ESP at the level of the state equation implies the existence of a state functional $H^F: {\cal Z}^{\mathbb{Z}_{-}} \longrightarrow {\cal X}$. Define ${\cal X}':=H^F({\cal Z}^{\mathbb{Z}_{-}})$. The relation \eqref{first for lemma} implies that $F$ restricts to a map $F: {\cal X}' \times {\cal Z} \longrightarrow {\cal X}' $ because for any $\mathbf{x}'=H^F({\bf z}') \in {\cal X}'  $,  with ${\bf z}' \in {\cal Z}^{\mathbb{Z}_{-}}  $ and any $\widetilde{ {\bf z} } \in {\cal Z} $,
\begin{equation*}
F(\mathbf{x}', \widetilde{ {\bf z} } )=H^F({\bf z}'\widetilde{ {\bf z} }) \in {\cal X}'.
\end{equation*}
The restricted state map obviously also has the ESP at the state level and has as associated functional the map with restricted codomain $H ^F:{\cal Z}^{\mathbb{Z}_{-}}  \longrightarrow {\cal X}' $, which proves that $({\cal X}',F,h ) $ is strongly reachable.
}

{Consider now a system $({\cal X},F,h ) $ such that  $F$ is invertible and $h$ is injective when restricted to ${\cal X}':=H^F({\cal Z}^{\mathbb{Z}_{-}})$. Let ${\bf z} _1, {\bf z} _2 \in {\cal Z} ^{\mathbb{Z}_{-}} $ be such that $H ^F_h({\bf z} _1\widetilde{{\bf z} })=H ^F_h({\bf z} _2\widetilde{{\bf z} })$ for any $\widetilde{{\bf z}} \in {\cal Z}^T$ and any $T \in \mathbb{N}$. The injectivity of $h$ implies that $H ^F({\bf z} _1\widetilde{{\bf z} })=H ^F({\bf z} _2\widetilde{{\bf z} })$.  Since the converse of \eqref{second for lemma} holds by the invertibility of $F$, we have  then that $H ^F({\bf z} _1)=H ^F({\bf z} _2)$ and we can hence conclude that the system does not have indistinguishable states and it is hence observable. \quad $\blacktriangledown$
}

\medskip

\noindent\textbf{Proof of Proposition \ref{proposition existence of morphism}.\ \ } Using the hypothesis on the strong reachability of $({\cal X},F,h ) $, we know that for any $\mathbf{x} \in {\cal X} $ there exists ${\bf z} \in {\cal Z}^{\mathbb{Z}_{-}}  $  such that $H ^F({\bf z})= \mathbf{x} $. Define:
\begin{equation*}
\begin{array}{rccc}
f: & {\cal X} & \longrightarrow &\overline{{\cal X}}\\
 &\mathbf{x}=H ^F({\bf z}) &\longmapsto & f(\mathbf{x}):=H^{\overline{F}}({\bf z}).
\end{array}
\end{equation*}
We now show that this map is well-defined and that it is the unique system morphism in the statement of the proposition. 

\begin{description}[leftmargin=1.5em]
\item [(i)] {\bf $f $ is well-defined:} given $\mathbf{x} \in {\cal X} $, let ${\bf z}_1, {\bf z}_2 \in {\cal Z}^{\mathbb{Z}_{-}} $ be such that $\mathbf{x}= H^F({\bf z} _1)= H^F({\bf z} _2) $. We now show that $H^{\overline{F}}({\bf z} _1)= H^{ \overline{F}}({\bf z} _2) $, necessarily. By contradiction, suppose that $\mathbf{x} _1:=H^{\overline{F}}({\bf z} _1)$,  $\mathbf{x} _2:= H^{ \overline{F}}({\bf z} _2) $, and that $\mathbf{x}_1\neq  {\mathbf{x} _2} $. As by hypothesis $(\overline{{\cal X}},\overline{F},\overline{h} ) $ is observable, there exists $\widetilde{{\bf z}} \in {\cal Z}^T  $, for some $T \in \mathbb{N} $, such that 
\begin{equation}
\label{out of observability}
H^{\overline{F}}_{\overline{h}}({\bf z} _1 \widetilde{{\bf z}})\neq H^{ \overline{F}}_{\overline{h}}({\bf z} _2 \widetilde{{\bf z}}). 
\end{equation}
However, the equality $H^F({\bf z} _1)= H^F({\bf z} _2) $ and \eqref{second for lemma} in Lemma \ref{preliminary for canonical} imply that $H^F({\bf z} _1 \widetilde{{\bf z}})= H^F({\bf z} _2 \widetilde{{\bf z}}) $ and hence $H^F_h({\bf z} _1 \widetilde{{\bf z}})= H^F_h({\bf z} _2 \widetilde{{\bf z}}) $. The hypothesis $H ^F_h=H ^{ \overline{F}}_{ \overline{h}} $ implies that $H ^{ \overline{F}}_{ \overline{h}} ({\bf z} _1 \widetilde{{\bf z}})= H ^{ \overline{F}}_{ \overline{h}}({\bf z} _2 \widetilde{{\bf z}}) $ which contradicts \eqref{out of observability}.
\item [(ii)] {\bf $f $ is system equivariant:} Let $\mathbf{x} \in {\cal X}  $, $ {\bf z} \in  {\cal Z}^{\mathbb{Z}_{-}}$, and $\widetilde{{\bf z}} \in {\cal Z} $, be such that $\mathbf{x}=H^F({\bf z}) $. Then, by \eqref{first for lemma} in Lemma \ref{preliminary for canonical} we have that
\begin{equation*}
f(F(\mathbf{x} , \widetilde{{\bf z}}))=f(F(H^F({\bf z}), \widetilde{{\bf z}}))=f(H^F({\bf z}\widetilde{{\bf z}}))=H^{ \overline{F}}({\bf z}\widetilde{{\bf z}})=\overline{F}(H^{\overline{F}}({\bf z}), \widetilde{{\bf z}})=\overline{F}(f (\mathbf{x}), \widetilde{{\bf z}}),
\end{equation*}
as required.
\item [(iii)] {\bf $f $ is readout invariant:} using the same elements as in the previous point:
\begin{equation*}
h(\mathbf{x})=h(H ^F({\bf z}))=H ^F_h({\bf z})=H^{\overline{F}}_{\overline{h}}({\bf z}) =\overline{h}(H ^{\overline{F}}({\bf z}))= \overline{h}(f (\mathbf{x})),
\end{equation*}
as required.
\item [(iv)] {\bf $f $ is unique:} Let $\overline{f}:\mathcal{X} \longrightarrow \overline{\mathcal{X}}$ be another system morphism. Let $\mathbf{x}=H^F({\bf z}) \in {\cal X} $  arbitrary. We first show that the sequence $\left(\overline{f}(H ^F(T_{-t}({\bf z}))), {\bf z} _t\right)_{t \in \mathbb{Z}_{-} } \in (\overline{{\cal X}} \times {\cal Z})^{\mathbb{Z}_{-}}$ is a solution of the system associated to $\overline{F} $. Indeed, for any $t \in \mathbb{Z}_{-} $, and by \eqref{first for lemma} and the system equivariance of $\overline{f} $:
\begin{equation*}
\overline{f}(H ^F(T_{-t}({\bf z})))=\overline{f}(F(U^F(T_{-t}({\bf z}))_{-1}, {\bf z} _t))=
\overline{f}(F(U^F(T_{-(t-1)}({\bf z}))_{0}, {\bf z} _t))=
\overline{F}(\overline{f}(H^F(T_{-(t-1)}({\bf z}))), {\bf z} _t),
\end{equation*}
as required. Now, since $\left((H ^{\overline{F}}(T_{-t}({\bf z}))), {\bf z} _t\right)_{t \in \mathbb{Z}_{-} } \in (\overline{{\cal X}} \times {\cal Z})^{\mathbb{Z}_{-}}$ is also a solution for the system associated to $\overline{F} $ that, by hypothesis, has the echo state property, we necessarily have that:
\begin{equation*}
\overline{f}(\mathbf{x})=\overline{f}(H^F({\bf z}))=H^{\overline{F}}({\bf z})= f(\mathbf{x}),
\end{equation*}
which proves the uniqueness of the morphism $f$. \quad $\blacksquare$
\end{description}

\begin{corollary}
\label{uniqueness corollary}
If the two systems $({\cal X},F,h ) $ and $(\overline{{\cal X}},\overline{F},\overline{h} ) $ in the statement of Proposition \ref{proposition existence of morphism} are canonical then they are necessarily system isomorphic.
\end{corollary}

\noindent\textbf{Proof.\ \ } By Proposition \ref{proposition existence of morphism}, the maps $f : {\cal X} \longrightarrow \overline{\mathcal{X}}$ and $\overline{f}: \overline{{\cal X}}  \longrightarrow {\cal X} $ defined by $f(\mathbf{x}):= H^{\overline{F}}({\bf z}) $, with $\mathbf{x} = H ^F({\bf z} )$, and $\overline{f}(\overline{\mathbf{x}}):=H ^F(\overline{{\bf z} }) $, with $\overline{\mathbf{x}}= H^{\overline{F}}(\overline{{\bf z}})$, for ${\bf z}, \overline{{\bf z}} \in \mathcal{Z}^{\mathbb{Z}_-}$, are well-defined system morphisms. Then, for any $\mathbf{x} = H ^F({\bf z} ) \in {\cal X}$ and $\overline{\mathbf{x}}= H^{\overline{F}}(\overline{{\bf z}})\in \overline{{\cal X}}$ we can verify that
\begin{equation*}
\overline{f} \circ  f (\mathbf{x}) = \overline{f}\left( H^{\overline{F}}({{\bf z}})\right)=H ^F({\bf z})= \mathbf{x}, \quad \mbox{and} \quad
f \circ \overline{f}(\overline{\mathbf{x}})=f \left(H ^F(\overline{{\bf z}})\right)=H^{\overline{F}}(\overline{{\bf z}})= \overline{\mathbf{x}},
\end{equation*}
which shows that $\overline{f}= f ^{-1} $ and $f= \overline{f} ^{-1} $, as required. \quad $\blacksquare$

\section{Canonical Realization Theorems}
\label{Canonical Realization Theorems}

In this section we propose two results in connection with the state-space system realization of input/output systems. The first result shows that any causal and time-invariant input/output system with discrete semi-infinite inputs admits a canonical state-space realization that is unique up to system isomorphisms. As we shall see later on in the examples in Section \ref{Realization and canonicalization of linear filters}, there is no guarantee that this realization takes place in a finite dimensional space. In a second result, we show that given any state-space system that satisfies the echo state property, we can always associate to it a canonical state-space realization (also unique up to system isomorphisms) that generates the same input/output system. This new canonical system is obtained from the original one by a procedure that we will generically call {\bfi  reduction} and is defined on a new state space  whose dimension (whenever that term is well-defined) is equal or smaller.

Apart from the causality and time-invariance, there is another dynamical feature that is needed to ensure the existence of these canonical realizations, namely, the {\bfi  input forgetting property} (see \cite{jaeger2001}).

\begin{definition}
\label{definition input forgetting property}
Let ${\cal Z} $ be a set, $({\cal Y}, d)$  a metric space, and let $U: {\cal Z}^{\mathbb{Z}_{-}} \longrightarrow {\cal Y}^{\mathbb{Z}_{-}} $ be a causal and time-invariant filter. We say that $U$ has the {\bf input forgetting property} whenever for any ${\bf u}, {\bf v} \in {\cal Z}^{\mathbb{Z}_{-}}$ and any ${\bf z} \in {\cal Z}^{\mathbb{N}^+}$:
\begin{equation}
\label{ifp}
\lim_{t \to \infty} d(H_U({\bf u \widetilde{z}}_t) , H_U({\bf v\widetilde{z}}_t)) = 0,
\end{equation}
where ${\bf \widetilde{z} _t} := \left({\bf z} _1, \ldots, {\bf z} _t\right)\in {\cal Z} ^t$, $t \in \mathbb{N}^+ $ {and $H _U: {\cal Z}^{\mathbb{Z}_{-}} \longrightarrow {\cal Y}  $ is the funcional associated to $U$ and defined by $H _U({\bf z})=U ({\bf z})_0 $.}
\end{definition}

This property is also referred to in the literature as the {\bfi unique steady-state property} (see \cite{Boyd1985}) and is usually obtained as a consequence of various continuity properties like the {\bfi  fading memory property} (see, for instance, \cite[Theorem 24]{RC9} and the definition later on in Section \ref{Realization and canonicalization of linear filters}).

\begin{theorem}[Canonical realization of input/output systems]
\label{Canonical realization of input/output systems}
Let ${\cal Z} $ be a set, $({\cal Y}, d)$  a metric space, and let $U: {\cal Z}^{\mathbb{Z}_{-}} \longrightarrow {\cal Y}^{\mathbb{Z}_{-}} $  be a causal and time-invariant input/output system that has the input forgetting property. Then, there exists a canonical state-space system $({\cal X}, F, h) $ such that $U=U^F_h $. This canonical realization of $U$ is unique up to system isomorphisms.
\end{theorem}

\noindent\textbf{Proof.\ \ } We start by defining the so-called {\bfi  Nerode equivalence relation} in ${\cal Z}^{\mathbb{Z}_{-}} $ with respect to the functional $H _U: {\cal Z} ^{\mathbb{Z}_{-}} \longrightarrow {\cal Y} $ determined by $U$ via the assignment $H _U({\bf z})= U ({\bf z}) _0$. We say that two elements ${\bf z} _1, {\bf z} _2 \in {\cal Z}^{\mathbb{Z}_{-}} $ are Nerode equivalent and write ${\bf z} _1\sim_{I} {\bf z} _2$, whenever $H _U({\bf z} _1\widetilde{{\bf z}})=H _U({\bf z} _2\widetilde{{\bf z}})$, for all $\widetilde{{\bf z}} \in {\cal Z} ^T  $ and all $T \in \mathbb{N} $. Define ${\cal X}:= {\cal Z}^{\mathbb{Z}_{-}}/\sim_{I} $, where the right-hand side of this equality stands for the set of equivalence classes in ${\cal Z}^{\mathbb{Z}_{-}} $ determined by the equivalence relation $\sim_{I} $, and denote by $[{\bf z}] \in {\cal X} $ the class that contains the element ${\bf z} \in {\cal Z}^{\mathbb{Z}_{-}} $.

Define now the system $({\cal X}, F,h) $, with $F: {\cal X} \times {\cal Z} \longrightarrow{\cal X} $ and  $h : {\cal X} \longrightarrow {\cal Y} $ given by 
\begin{equation}
\label{realizing state system}
F([\mathbf{z}], \widetilde{{\bf z}}):= [{\bf z} \widetilde{{\bf z}}] \quad \mbox{and} \quad
h([{\bf z}]):= H _U ({\bf z}). 
\end{equation}
We now show that this system is well-defined, it has the echo state property, and that it is a canonical realization of $U $. If that is the case, the uniqueness up to system isomorphisms follows from Corollary \ref{uniqueness corollary}. We proceed point by point:
\begin{description}[leftmargin=1.5em]
\item [(i)] {\bf $({\cal X}, F, h)$ is well-defined:} First of all, $F: {\cal X} \times {\cal Z} \longrightarrow{\cal X} $ is well-defined because if ${\bf z} _1, {\bf z}_2 \in {\cal Z} ^{\mathbb{Z}_{-}}$  are such that  ${\bf z} _1\sim_{I} {\bf z} _2 $ then, by definition, 
\begin{equation}
\label{nerode for later}
H _U({\bf z} _1\widetilde{{\bf z}})=H _U({\bf z} _2\widetilde{{\bf z}}), \quad \mbox{for all $\widetilde{{\bf z}} \in {\cal Z} ^T  $ and all $T \in \mathbb{N} $. }
\end{equation}
In particular, for any $\widehat{{\bf z}} \in {\cal Z} $, we have that $F([{\bf z} _1], \widehat{{\bf z}})=   F([{\bf z} _2], \widehat{{\bf z}})$ because $[{\bf z} _1\widehat{{\bf z}}]=[{\bf z} _2\widehat{{\bf z}}] $, as \eqref{nerode for later} also implies that $H _U({\bf z} _1\widehat{{\bf z}}\widetilde{{\bf z}})=H _U({\bf z} _2\widehat{{\bf z}}\widetilde{{\bf z}}) $ for all $\widetilde{{\bf z}} \in {\cal Z} ^T  $ and all $T \in \mathbb{N} $. The map $h : {\cal X} \longrightarrow {\cal Y}$ is also well-defined because if we consider ${\bf z} _1, {\bf z}_2 \in {\cal Z} ^{\mathbb{Z}_{-}}$ that, as above, ${\bf z} _1\sim_{I} {\bf z} _2 $, the equality \eqref{nerode for later} implies, in particular, that $H _U({\bf z} _1)=H _U({\bf z} _2)$ and hence $h([{\bf z} _1])=H _U({\bf z} _1)=H _U({\bf z} _2)=h([{\bf z} _2]) $.
\item [(ii)] {\bf The system $({\cal X}, F, h)$ has the echo state property:} Given ${\bf z} \in {\cal Z} ^{\mathbb{Z}_{-}}  $, we first show that the sequence $ \left([T_{-t}({\bf z})], {\bf z} _t\right)_{t \in \mathbb{Z}_{-}}\in ({\cal X} \times {\cal Z})^{\mathbb{Z}_{-}}$ is a solution of the state system $({\cal X}, F)$. This is so because, for any $t \in \mathbb{Z}_{-} $, we have
\begin{equation*}
F \left([T_{-(t-1)}({\bf z})], {\bf z} _t\right)= \left[T_{-(t-1)}({\bf z}) {\bf z} _t  \right]= [T_{-t}({\bf z})].
\end{equation*}
We now show that this solution is unique. Suppose that $(\mathbf{x} _t, {\bf z} _t)_{t \in \mathbb{Z}_{-}} \in ({\cal X} \times {\cal Z})^{\mathbb{Z}_{-}} $ is also a solution for $({\cal X}, F)$ with respect to the same input sequence. Since the quotient map $ {\cal Z}^{\mathbb{Z}_{-}} \longrightarrow {\cal Z}^{\mathbb{Z}_{-}}/\sim_{I} $ is surjective, for any $t \in \mathbb{Z}_{-} $  there exists an element $\overline{{\bf z}} _t \in {\cal Z}^{\mathbb{Z}_{-}}  $  such that $\mathbf{x} _t= [\overline{{\bf z}} _t ]$. The solution condition on $(\mathbf{x} _t, {\bf z} _t)_{t \in \mathbb{Z}_{-}} $  implies that, also for any $t \in \mathbb{Z}_{-} $, $[\overline{{\bf z}} _t ]=F([\overline{{\bf z}} _{t-1} ], {\bf z} _t)=[\overline{{\bf z}} _{t-1}  {\bf z} _t] $ and hence $H _U(\overline{{\bf z}} _t \widetilde{{\bf z}})=H _U(\overline{{\bf z}} _{t-1}  {\bf z} _t\widetilde{{\bf z}}) $, for all $\widetilde{{\bf z}} \in {\cal Z} ^T  $ and all $T \in \mathbb{N} $. If we use recursively this identity, we can show that
\begin{equation*}
H _U(\overline{{\bf z}} _t \widetilde{{\bf z}})=H _U(\overline{{\bf z}} _{t-1}  {\bf z} _t\widetilde{{\bf z}})=H _U(\overline{{\bf z}} _{t-2}  {\bf z} _{t-1} {\bf z} _t\widetilde{{\bf z}})= \cdots =H _U(\overline{{\bf z}} _{t-\tau} {{\bf z}} _{t-(\tau-1)} \cdots {\bf z} _{t-1}  {\bf z} _t\widetilde{{\bf z}}), 
\end{equation*}
for all $\widetilde{{\bf z}} \in {\cal Z} ^T  $ and all $\tau,T \in \mathbb{N} $. 
These equalities imply that for any $t \in \Bbb Z_- $  and $\tau\in \mathbb{N} $:
\begin{equation*}
d\left(H _U(\overline{{\bf z}} _t \widetilde{{\bf z}}),H_U(T_{-t}({\bf z})\widetilde{{\bf z}})\right)=
d\left(H _U(\overline{{\bf z}} _{t-\tau} {{\bf z}} _{t-(\tau-1)} \cdots {\bf z} _{t-1}  {\bf z} _t\widetilde{{\bf z}}),H_U(T_{-t}({\bf z})\widetilde{{\bf z}})\right).
\end{equation*}
Now, since by hypothesis $U $ satisfies the input forgetting property, we can take a limit on $\tau  $  on the right-hand side of this equality and conclude that
\begin{equation*}
d\left(H _U(\overline{{\bf z}} _t \widetilde{{\bf z}}),H_U(T_{-t}({\bf z})\widetilde{{\bf z}})\right)=\lim_{\tau \rightarrow \infty}
d\left(H _U(\overline{{\bf z}} _{t-\tau} {{\bf z}} _{t-(\tau-1)} \cdots {\bf z} _{t-1}  {\bf z} _t\widetilde{{\bf z}}),H_U(T_{-t}({\bf z})\widetilde{{\bf z}})\right)=0,
\end{equation*}
which implies that $H _U(\overline{{\bf z}} _t \widetilde{{\bf z}})=H_U(T_{-t}({\bf z})\widetilde{{\bf z}})$ and hence that $\mathbf{x} _t= [\overline{{\bf z}} _t ]=[T_{-t}({\bf z})]$, as required.
\item [(iii)] {\bf $({\cal X}, F, h)$ is a state-space realization of $U$:} Since in the previous point we proved that $({\cal X}, F, h)$ has the echo state property, we can  associate to it a system filter $U ^F _h: {\cal Z} ^{\mathbb{Z}_{-}} \longrightarrow {\cal Y}^{\mathbb{Z}_{-} }$. We also showed that for any input ${\bf z} \in {\cal Z} ^{\mathbb{Z}_{-}} $ the sequence $ \left([T_{-t}({\bf z})], {\bf z} _t\right)_{t \in \mathbb{Z}_{-}}\in ({\cal X} \times {\cal Z})^{\mathbb{Z}_{-}}$ is the unique solution of the state system $({\cal X}, F)$ which proves that the  state filter $U ^F: {\cal Z} ^{\mathbb{Z}_{-}} \longrightarrow {\cal X}^{\mathbb{Z}_{-} }$ is given by
\begin{equation}
\label{state canonical filter}
U ^F({\bf z})_t=[T_{-t}({\bf z})].
\end{equation}
Consequently, for any $t \in \mathbb{Z}_{-} $, we have that
\begin{equation}
\label{state canonical system}
U ^F _h( {\bf z})_t=h \left([T_{-t}({\bf z})]\right)=H _U(T_{-t}({\bf z}))=U ({\bf z})_t,
\end{equation}
which implies that $U ^F _h=U $.
\item [(iv)] {\bf $({\cal X}, F, h)$ is canonical:} Since for any ${\bf z}\in {\cal Z}^{\mathbb{Z}_{-}} $ the equality \eqref{state canonical filter} guarantees that $H^F( {\bf z})= [ {\bf z}] $, we can immediately conclude that $({\cal X}, F, h)$ is strongly reachable. Let now $\mathbf{x} _1= H ^F({\bf z} _1)=[{\bf z} _1]$ and $\mathbf{x} _2= H ^F({\bf z} _2)=[{\bf z} _2]$ be two indistinguishable states, that is,  for any $\widetilde{{\bf z}} \in \Bbb Z ^T$ and any $T \in \mathbb{N} $, we have that $H ^F_h({\bf z} _1\widetilde{{\bf z}})=H ^F_h({\bf z} _2\widetilde{{\bf z}})$. The equality \eqref{state canonical system} evaluated at $t =0  $ implies that in that case $H_U({\bf z} _1\widetilde{{\bf z}})=H _U({\bf z} _2\widetilde{{\bf z}})$, necessarily, and hence we can conclude that $[{\bf z} _1]=[{\bf z} _2  ] $, which is equivalent to $\mathbf{x} _1= \mathbf{x}_2 $, as required. \quad $\blacksquare$
\end{description}

\begin{remark}
\normalfont
It is easy to see that Theorem \ref{Canonical realization of input/output systems} remains valid when the spaces $ {\cal Z}^{\mathbb{Z}_{-}} $  and $ {\cal Y}^{\mathbb{Z}_{-}}  $ are replaced by time-invariant subsets $\mathcal{V}_{{\cal Z}}\subset  {\cal Z}^{\mathbb{Z}_{-}} $  and $\mathcal{V}_{{\cal Y}}\subset {\cal Y}^{\mathbb{Z}_{-}}  $, respectively, that additionally are also invariant with respect to the concatenation with finite sequences that was used in the definition of the Nerode equivalence relation. The time invariance is defined by the property $T _\tau(\mathcal{V}_{{\cal Z}}) \subset \mathcal{V}_{{\cal Z}}  $ and $T _\tau(\mathcal{V}_{{\cal Y}}) \subset \mathcal{V}_{{\cal Y}}  $, for any $\tau \in \mathbb{N} $.
\end{remark}

The canonicalization theorem that we just proved provides a canonical state-space realization for any input-forgetting, causal, and time-invariant filter by using as state-space the set of equivalence classes in the space of semi-infinite input sequences with respect to the Nerode equivalence. If that filter happens to be already given in a state-space form, we shall show in the next theorem that a canonical realization can be constructed for it by {\bfi  reducing} the given state-space. 

The reduction procedure  that we propose next is reminiscent of the {\it optimal reduction} method introduced in \cite{optimal:mm, reduction:optimal:cras} in the context of symmetric Hamiltonian systems and consists in two steps. First, given a (generically non-canonical) state-space system $({\cal X}, F, h)$ with ${\cal Z}$ and $ {\cal Y} $ as input and output spaces, respectively, and  that satisfies the echo state property, we restrict the state equation to the subset ${\cal X}_R \subset {\cal X} $ of {\bfi  reachable states} defined by
\begin{equation}
\label{set of reachable states}
{\cal X} _R:= \left\{\mathbf{x} \in {\cal X}\mid \mathbf{x}= H ^F({\bf z})\  \mbox{for some} \  {\bf z} \in {\cal Z}^{\mathbb{Z}_{-}}\right\}.
\end{equation}
Note that ${\cal X}_R $ is the state subspace already introduced in Corollary \ref{ESP implies exists reachable set}.

In a second step, we can define in ${\cal X}_R $ the Nerode equivalence relation $\sim_{S} $ that in the previous theorem was formulated in the space of semi-infinite input sequences. More explicitly, given $\mathbf{x} _1= H ^F({\bf z} _1), \mathbf{x} _2= H ^F({\bf z} _2) \in {\cal X} _R $, for some ${\bf z} _1, {\bf z} _2 \in  {\cal Z} ^{\mathbb{Z}_{-}}  $, we say that these two states are {\bfi  Nerode equivalent} and, as before, we denote
\begin{equation}
\label{definition nerode in states}
\mathbf{x} _1\sim_{S} \mathbf{x} _2 \  \mbox{whenever $H _h ^F({\bf z} _1\widetilde{{\bf z}})=H  _h ^F({\bf z} _2\widetilde{{\bf z}})$, for all $\widetilde{{\bf z}} \in {\cal Z} ^T  $ and all $T \in \mathbb{N} $.  }
\end{equation}
Notice that this definition of Nerode equivalent states is equivalent to the so-called indistinguishable states which is introduced in part {\bf (ii)} of Definition \ref{Reachable and Observable}.

The symbol $[\mathbf{x}] \in {\cal X}_R/\sim_{S} $ denotes the equivalence class that contains the element $\mathbf{x} \in {\cal X}_R $. We emphasize that this relation is well-defined since it does not depend on the elements ${\bf z} _1, {\bf z} _2 \in  {\cal Z} ^{\mathbb{Z}_{-}}  $ used to define $\mathbf{x} _1 $ and $\mathbf{x} _2 $ because of \eqref{second for lemma} in Lemma \ref{preliminary for canonical}.

In the next theorem will show that $({\cal X}, F, h)$ naturally projects to a system on the quotient ${\cal X}_R/\sim_{S} $ that has the echo state property if $({\cal X}, F, h)$ is input-forgetting and, more importantly, is canonical.

\begin{theorem}[Canonicalization by reduction]
\label{Canonicalization by reduction}
Let ${\cal Z} $ be a set, $({\cal Y}, d)$  a metric space, and let $({\cal X}, F, h)$ be a state-space system that has ${\cal Z} $ and ${\cal Y}$ as input and output spaces, respectively. Suppose that $({\cal X}, F) $ has the echo state property and that the state-space  filter $U^F_h: {\cal Z} ^{\mathbb{Z}_{-}} \longrightarrow {\cal Y}^{\mathbb{Z}_{-}}  $ has the input forgetting property. Let ${\cal X}_R \subset {\cal X} $ be the set of reachable states defined in \eqref{set of reachable states} and $\overline{{\cal X}}:={\cal X}_R/\sim_{S} $ the quotient set with respect to the Nerode equivalence relation $\sim_{S} $ defined in \eqref{definition nerode in states}. 

The state-space system $({\cal X}, F, h)$ drops to another system $(\overline{{\cal X}}, \overline{F}, \overline{h})$ with the same input and output spaces, with states in the quotient space $\overline{{\cal X}} $, and maps $\overline{F}: \overline{{\cal X}} \times {\cal Z} \longrightarrow \overline{{\cal X}} $ and $\overline{h}: \overline{{\cal X}} \longrightarrow{\cal Y}  $ defined by:
\begin{equation}
\label{definition reduced system}
\left\{
\begin{array}{lll}
\overline{F}([\mathbf{x}], {\bf z})&:= &[F(\mathbf{x}, {\bf z})],\\
\overline{h}([\mathbf{x}]) &:= & h( \mathbf{x}).
\end{array}
\right.
\end{equation}
The state-space system $(\overline{{\cal X}}, \overline{F}, \overline{h})$ has the echo state property and it is a canonical realization of $U^F_h$. We refer to $(\overline{{\cal X}}, \overline{F}, \overline{h})$ as the {\bfi  canonical reduced realization} of  $({\cal X}, F, h)$.
\end{theorem}

\noindent\textbf{Proof.\ \ }We first show that the reduced state and readout maps $\overline{F} $ and $\overline{h} $ in \eqref{definition reduced system} are well-defined. Concerning $\overline{F } $, we show first that the restriction of $F $ to ${\cal X} _R \times {\cal Z}  $  maps into ${\cal X} _R $. Indeed, let $\mathbf{x}  \in {\cal X} _R $ arbitrary and let ${\bf z} \in {\cal Z}^{\mathbb{Z}_{-}} $ be such that  $\mathbf{x} =H ^F({\bf z})$. Then, for any $\widetilde{{\bf z}} \in {\cal Z} $, by \eqref{first for lemma} in Lemma \ref{preliminary for canonical}, we have that
\begin{equation*}
F(\mathbf{x}, \widetilde{{\bf z}})=F(H ^F({\bf z}), \widetilde{{\bf z}})=H^F({\bf z}\widetilde{{\bf z}}) \in {\cal X} _R.
\end{equation*}  
This guarantees that $F: {\cal X} \times {\cal Z} \longrightarrow {\cal X} $ restricts to a map $F _R: {\cal X}_R \times {\cal Z} \longrightarrow{\cal X} _R $ that we now show drops to $\overline{F}: \overline{{\cal X}} \times {\cal Z} \longrightarrow \overline{{\cal X}} $ by proving that if $ \mathbf{x} _1, \mathbf{x} _2 \in {\cal X} _R $ are such that $ \mathbf{x} _1\sim_{S} \mathbf{x} _2 $, then $F _R(\mathbf{x}_1, {\bf z})\sim_{S} F _R(\mathbf{x}_2, {\bf z}) $, for all ${\bf z} \in {\cal Z} $. Indeed, if $ \mathbf{x} _1\sim_{S} \mathbf{x} _2 $, by definition \eqref{definition nerode in states}, $H _h ^F({\bf z} _1\widetilde{{\bf z}})=H  _h ^F({\bf z} _2\widetilde{{\bf z}})$, for all $\widetilde{{\bf z}} \in {\cal Z} ^T  $ and all $T \in \mathbb{N} $, where $\mathbf{x} _1= H ^F({\bf z} _1), \mathbf{x} _2= H ^F({\bf z} _2) \in {\cal X} _R $, for some ${\bf z} _1, {\bf z} _2 \in  {\cal Z} ^{\mathbb{Z}_{-}}  $. Now, by \eqref{first for lemma} and for all ${\bf z} \in {\cal Z} $, $F _R(\mathbf{x} _1, {\bf z})= H^F({\bf z}_1 {\bf z}) $, $F _R(\mathbf{x} _2, {\bf z})= H^F({\bf z}_2 {\bf z}) $ and since by  \eqref{definition nerode in states} $H _h ^F({\bf z} _1{\bf z}\widetilde{{\bf z}})=H  _h ^F({\bf z} _2 {\bf z}\widetilde{{\bf z}})$, for all $\widetilde{{\bf z}} \in {\cal Z} ^T  $ and all $T \in \mathbb{N} $, we can conclude that $F _R(\mathbf{x}_1, {\bf z})\sim_{S} F _R(\mathbf{x}_2, {\bf z}) $, as required. In order to show that $\overline{h} $ is well-defined, consider first the restriction $h _R:=h\mid _{{\cal X}_R} : {\cal X} _R \longrightarrow {\cal Y}$ as well as two elements $ \mathbf{x} _1, \mathbf{x} _2 \in {\cal X} _R $ as above such that $ \mathbf{x} _1\sim_{S} \mathbf{x} _2 $. Taking now for $\widetilde{\bf z} $ the empty sequence in the definition of the equivalence relation $\sim_{S} $, we have that:
\begin{equation*}
h_R(\mathbf{x} _1)=h \left(H ^F({\bf z} _1)\right)=H^F_{h}({\bf z} _1)=H^F_{h}({\bf z} _2)=h \left(H ^F({\bf z} _2)\right)=h_R(\mathbf{x} _2),
\end{equation*}
which proves that $h_R $  drops to the map $\overline{h} $ in the statement and it is hence well-defined.

We now show that the reduced system $(\overline{{\cal X}}, \overline{F}, \overline{h})$ has the echo state property by following a scheme similar to part {\bf (ii)} in the proof of Theorem \ref{Canonical realization of input/output systems}. First of all, it is easy to see that if  $({\bf x} _t,{\bf z} _t)_{t \in \mathbb{Z}_{-}}$ is the unique solution of the system $({\cal X},F)$ (that by hypothesis satisfies the echo state property) associated to ${\bf z}:=({\bf z} _t)_{t \in \mathbb{Z}_{-}}$, then $([{\bf x} _t],{\bf z} _t)_{t \in \mathbb{Z}_{-}}$ is a  solution of the system $(\overline{{\cal X}},\overline{F})$ associated to ${\bf z}$. We now show that that solution is unique. Suppose that  $([\overline{{\bf x} }_t],{\bf z} _t)_{t \in \mathbb{Z}_{-}}$ is another solution of $(\overline{{\cal X}},\overline{F})$ for the same input ${\bf z}$. For any $t \in \mathbb{Z}_{-} $,  let $\overline{{\bf z}} _t \in {\cal Z}^{\mathbb{Z}_{-}}  $ be such that $\overline{\mathbf{x}} _t= H^F(\overline{{\bf z}} _t)$. The solution condition on $([\overline{\mathbf{x}} _t], {\bf z} _t)_{t \in \mathbb{Z}_{-}} $  implies that, also for any $t \in \mathbb{Z}_{-} $, 
\begin{equation*}
[H^F(\overline{{\bf z}} _t)]=[\overline{\mathbf{x}} _t]= \overline{F}([\overline{\mathbf{x}}_{t-1}], {\bf z} _t)=[ {F}(\overline{\mathbf{x}}_{t-1}, {\bf z} _t)]= [H^F(\overline{{\bf z}} _{t-1}{\bf z} _t)],
\end{equation*}
which by \eqref{definition nerode in states},  implies that for all $\widetilde{{\bf z}} \in {\cal Z} ^T  $ and all $T \in \mathbb{N} $ one has $H _h ^F(\overline{{\bf z}} _{t}\widetilde{{\bf z}})=H  _h ^F(\overline{{\bf z}} _{t-1}{\bf z} _t\widetilde{{\bf z}})$, necessarily. 
If we use recursively this identity, we can show that
\begin{equation*}
H _h^F(\overline{{\bf z}} _t \widetilde{{\bf z}})=H _h^F(\overline{{\bf z}} _{t-1}  {\bf z} _t\widetilde{{\bf z}})=H _h^F(\overline{{\bf z}} _{t-2}  {\bf z} _{t-1} {\bf z} _t\widetilde{{\bf z}})= \cdots =H _h^F(\overline{{\bf z}} _{t-\tau} {{\bf z}} _{t-(\tau-1)} \cdots {\bf z} _{t-1}  {\bf z} _t\widetilde{{\bf z}}), 
\end{equation*}
for all $\widetilde{{\bf z}} \in {\cal Z} ^T  $ and all $\tau,T \in \mathbb{N} $. 
These equalities imply that for any $t \in \Bbb Z_- $  and $\tau\in \mathbb{N} $:
\begin{equation*}
d\left(H _h^F(\overline{{\bf z}} _t \widetilde{{\bf z}}),H_h^F(T_{-t}({\bf z})\widetilde{{\bf z}})\right)=
d\left(H _h^F(\overline{{\bf z}} _{t-\tau} {{\bf z}} _{t-(\tau-1)} \cdots {\bf z} _{t-1}  {\bf z} _t\widetilde{{\bf z}}),H_h^F(T_{-t}({\bf z})\widetilde{{\bf z}})\right).
\end{equation*}
Now, since by hypothesis $H_h^F$ has the input forgetting property, we can take a limit on $\tau  $  on the right-hand side of this equality and conclude that
\begin{equation*}
d\left(H _h^F(\overline{{\bf z}} _t \widetilde{{\bf z}}),H_h^F(T_{-t}({\bf z})\widetilde{{\bf z}})\right)=\lim_{\tau \rightarrow \infty}
d\left(H _h^F(\overline{{\bf z}} _{t-\tau} {{\bf z}} _{t-(\tau-1)} \cdots {\bf z} _{t-1}  {\bf z} _t\widetilde{{\bf z}}),H_h^F(T_{-t}({\bf z})\widetilde{{\bf z}})\right)=0,
\end{equation*}
which implies that $H _h^F(\overline{{\bf z}} _t \widetilde{{\bf z}})=H_h^F(T_{-t}({\bf z})\widetilde{{\bf z}})$ and hence that $[\overline{\mathbf{x}} _t]= [H ^F(\overline{{\bf z}} _t )]=[H ^F(T_{-t}({\bf z}))]=[\mathbf{x} _t]$, as required. 

Finally, the fact that $([{\bf x} _t],{\bf z} _t)_{t \in \mathbb{Z}_{-}}$ is the unique  solution of  $(\overline{{\cal X}},\overline{F})$ associated to ${\bf z}$ when $({\bf x} _t,{\bf z} _t)_{t \in \mathbb{Z}_{-}}$ is the unique solution of  $({\cal X},F)$ amounts to the equality $H^{\overline{F}}_{\overline{h}} =H^F_h$. Consequently, $(\overline{{\cal X}}, \overline{F}, \overline{h}) $ is a realization for the filter associated to $({\cal X},F,h) $ and it is trivially canonical. \quad $\blacksquare$

\medskip

Since Theorems  \ref{Canonical realization of input/output systems} and \ref{Canonicalization by reduction}  produce two different canonical realizations of a given system and we know by Corollary \ref{uniqueness corollary} that those realizations are unique up to system isomorphisms, we can conclude the non-trivial statement that the two sets of classes $ {\cal Z}^{\mathbb{Z}_{-}}/\sim_{I} $ and ${\cal X}_R/\sim_{S} $ in the space of semi-infinite input sequences and on the space of reachable states, respectively, are isomorphic quotient spaces. We frame that result in the next corollary.

\begin{corollary}
Let ${\cal Z} $ be a set, $({\cal Y}, d)$  a metric space, and let $({\cal X}, F, h)$ be a state-space system that has ${\cal Z} $ and ${\cal Y}$ as input and output spaces, respectively. Suppose that $({\cal X}, F) $ has the echo state property and that the state-space system filter $U^F_h: {\cal Z} ^{\mathbb{Z}_{-}} \longrightarrow {\cal Y}^{\mathbb{Z}_{-}}  $ has the input forgetting property. Let ${\cal X}_R/\sim_{S} $ be the reduced state-space defined in \eqref{definition nerode in states} and let ${\cal Z}^{\mathbb{Z}_{-}}/\sim_{I} $ be the quotient space defined in the proof of Theorem \ref{Canonical realization of input/output systems}. These two quotient spaces are isomorphic. The isomorphism is implemented by the map:
\begin{equation*}
\begin{array}{cccc}
f: &{\cal Z}^{\mathbb{Z}_{-}}/\sim_{I} & \longrightarrow &{\cal X}_R/\sim_{S} \\
	&[ {\bf z}] &\longmapsto & [H ^F( {\bf z})].
\end{array}
\end{equation*}
\end{corollary}

\section{Realization and canonicalization of linear filters}
\label{Realization and canonicalization of linear filters}

In this section we study the realization and canonicalization problem for linear,  time-invariant, and causal filters that satisfy the so-called fading memory property. In order to explicitly define the input spaces and this property we first consider the {\bfi  supremum norm}
$\left\|\cdot \right\|_{\infty} $ in the space of semi-infinite sequences $\mathbb{R}^{\mathbb{Z}_{-}} $ in $\mathbb{R}$ defined by
\begin{equation}
\label{infinity norm}
\left\|{\bf z} \right\|_{\infty}:=\sup_{t \in \mathbb{Z}_{-}} \left\{|{ z} _t|\right\}, \quad \mbox{for any} \quad {\bf z}\in \mathbb{R}^{\mathbb{Z}_{-}}.
\end{equation}
Let $(\ell_{-} ^{\infty}(\mathbb{R}),  \left\|\cdot \right\|_ \infty) $ be the Banach space formed by the elements in $\mathbb{R}^{\mathbb{Z}_{-}} $ that have a finite supremum norm. We define now a {\bfi  weighting sequence} $w : \mathbb{N} \longrightarrow (0,1] $ as a a strictly decreasing sequence with zero limit such that $w _0=1 $. Given an element ${\bf z} \in \ell_{-} ^{\infty}(\mathbb{R}) $, we define its $w$-{\bfi weighted norm} $\left\|\cdot \right\|  _w $ by
\begin{equation*}
\left\|{\bf z} \right\|_{w}:=\sup_{t \in \mathbb{Z}_{-}} \left\{|{\bf z} _t|w _{-t}\right\}.
\end{equation*}

Consider now a linear,  time-invariant, and causal filter $U : \ell^{\infty}_-(\mathbb{R}) \longrightarrow \ell^{\infty}_-(\mathbb{R}) $. We say that the functional $H _U: \ell^{\infty}_-(\mathbb{R}) \longrightarrow \mathbb{R}$ associated to $U$ has the so-called {\bfi  fading memory property (FMP)} with respect to the weighting sequence  $w$ whenever for any $\epsilon>0 $, there exists $\delta(\epsilon)>0 $ such that if ${\bf z} \in \ell^{\infty}_-(\mathbb{R}) $ is such that $\left\|{\bf z}\right\|_w< \delta (\epsilon) $ then $|H _U({\bf z})|< \epsilon $, necessarily. 

The Convolution Theorem (see \cite[Theorem 5]{Boyd1985}) shows that $H _U $ has the FMP if and only if the filter $U$ has a convolution representation, that is, there exists an element $\boldsymbol{\Psi} \in \ell_{-}^{1}(\mathbb{R})$ such that 
\begin{equation*}
U({\bf z})_t=\sum_{j\in \mathbb{Z}_{-}} {\Psi} _j z_{t+j}=: \left(\boldsymbol{\Psi} \ast {\bf z}\right)_t, \quad \mbox{for any} \quad {\bf z} \in \ell^{\infty}_-(\mathbb{R}), \ t \in \mathbb{Z}_{-}.
\end{equation*}
In such case, it is easy to see that $U : \ell^{\infty}_-(\mathbb{R}) \longrightarrow \ell^{\infty}_-(\mathbb{R}) $ is a bounded linear operator and that its operator norm $\vertiii{U}_{\infty} $ satisfies that $\vertiii{U}_{\infty} \leq \left\|\boldsymbol{\Psi}\right\|_1=:\sum _{j \in \mathbb{Z}_{-}}|{\Psi} _j |< \infty$.

Additionally, as we already mentioned after Definition \ref{definition input forgetting property}, the FMP implies the input forgetting property that we used in the main results in Section \ref{Canonical Realization Theorems}. Since this fact is proved in the literature (see \cite[Theorem 6]{Boyd1985} and \cite[Theorem 6]{RC9}) exclusively for uniformly bounded inputs, we prove it separately in our situation in the following result that collects all the facts that we just mentioned. Before we proceed with the statement, we extend the definition of the time delay operator $T _\tau : {\cal Z}^{\mathbb{Z}_{-}}\longrightarrow {\cal Z}^{\mathbb{Z}_{-}} $ defined in \eqref{first time delay operator} for any $\tau \in \mathbb{N} $, to accommodate any $\tau \in \mathbb{Z} $ by setting, for any $\tau<0$:
\begin{equation}
\label{time delay negative}
T _{\tau} ({\bf z}):=({\bf z}, \underbrace{{\bf 0}, \ldots, {\bf 0}}_{\mbox{$-\tau$ times}}), \quad {\bf z} \in {\cal Z}^{\Bbb Z _-}.
\end{equation}

\begin{proposition}
\label{properties fmp filter linear}
Let $U : \ell^{\infty}_-(\mathbb{R}) \longrightarrow \ell^{\infty}_-(\mathbb{R}) $ be a linear, time-invariant, and causal filter such that $H _U: \ell^{\infty}_-(\mathbb{R}) \longrightarrow \mathbb{R} $ has the fading memory property with respect to a weighting sequence $w$. Then, there exists a unique element $\boldsymbol{\Psi} \in \ell_{-}^{1}(\mathbb{R}) $ such that  $U ({\bf z})= \boldsymbol{\Psi}\ast {\bf z} $ for any ${\bf z} \in \ell^{\infty}_-(\mathbb{R}) $. Moreover, $U $ is a bounded linear automorphism of $\ell^{\infty}_-(\mathbb{R}) $ such that  $\vertiii{U}_{\infty}\leq \left\|\boldsymbol{\Psi}\right\|_1< \infty $ and it has the input forgetting property.
\end{proposition}

\noindent\textbf{Proof.\ \ } In view of the references quoted above, it just remains to be shown that the element $\boldsymbol{\Psi} \in \ell_{-}^{1}(\mathbb{R}) $ that provides the convolution representation is unique and that $U$ has the input forgetting property.

The uniqueness of the sequence $\boldsymbol{\Psi} \in \ell_{-}^{1}(\mathbb{R}) $ is due to the fact that its components are uniquely determined by the impulse response of $U$, that is, for any $t \in \mathbb{Z}_{-}$
\begin{equation*}
\Psi_t=H _U(\mathbf{e}_t)= \left(\boldsymbol{\Psi} \ast \mathbf{e}_t\right)_0, \quad \mbox{where} \quad \mathbf{e}_t:=  (\ldots, 0,\underbrace{1}_{\mbox{$t$ entry}},0, \ldots,0 ) \in \ell^{\infty}_-(\mathbb{R}). 
\end{equation*}

We now show that $U$ has the input forgetting property. Let ${\bf u}, {\bf v} \in \ell^{\infty}_-(\mathbb{R})$, ${\bf z} \in \mathbb{R}^{\mathbb{N}^+}$, and denote ${\bf \widetilde{z} _t} := \left({ z} _1, \ldots, { z} _t\right)\in {\mathbb R} ^t$, for any $t \in \mathbb{N}^+ $. It is easy to see using \eqref{time delay negative} that
\begin{equation*}
\mathbf{u}{\bf \widetilde{z} _t}=T _{-t}(\mathbf{u})+ \left(\ldots,0,z _1, z _2, \ldots, z _t\right)\quad \mbox{and} \quad \mathbf{v}{\bf \widetilde{z} _t}=T _{-t}(\mathbf{v})+ \left(\ldots,0,z _1, z _2, \ldots, z _t\right).
\end{equation*}
We now use these equalities with the convolution representation of $U$ and the linearity of $T _{-t} $ and show that:
\begin{multline}
\label{implies fmp with l1}
\left|H _U(\mathbf{u}{\bf \widetilde{z} _t})-H _U(\mathbf{v}{\bf \widetilde{z} _t})\right|=\left| U(\mathbf{u}{\bf \widetilde{z} _t})_0-U(\mathbf{v}{\bf \widetilde{z} _t})_0\right|=\left|\boldsymbol{\Psi}\ast T _{-t}(\mathbf{u}-\mathbf{v})\right|=\left|\sum_{j=- \infty}^{-t}\Psi _j \left(u-v\right)_{j+t}\right|\\
\leq \sum_{j=- \infty}^{-t}\left|\Psi _j\right| \left\|\mathbf{u}-\mathbf{v}\right\| _{\infty}= \left( \left\|\boldsymbol{\Psi}\right\|_1-\sum_{j=-t+1}^0\left|\Psi _j\right|\right)\left\|\mathbf{u}-\mathbf{v}\right\| _{\infty}\xrightarrow[t \to \infty]{} 0. \quad \blacksquare
\end{multline}

\medskip

The proposition that we just proved shows, in particular, that FMP linear, causal, and time invariant filters satisfy the hypotheses of the Canonical Realization Theorem \ref{Canonical realization of input/output systems} and hence they always have a canonical state-space realization that, as we show later in Theorem \ref{Canonicalization of linear filters}, is linear even though the state space may be an infinite dimensional vector space. {We emphasize that the fading memory property plays a crucial role in the result that we just proved since, in general, the Convolution Theorem does not hold in its absence (see the counterexample in Section A3 of \cite{Boyd1985}).}

Before we proceed with that theorem, we first state a result that lists important properties of finite-dimensional linear state-space realizations that are needed in the sequel.

\begin{proposition}[Linear state-space realizations with semi-infinite inputs]
\label{Linear state-space realizations with semi-infinite inputs}
Let $N \in \mathbb{N} $, let $A \in \mathbb{M}_N  $ be a diagonalizable matrix, $\mathbf{C} \in \mathbb{R}^N  $, $W \in \mathbb{M}_{1,N} $, and consider the linear state-space system $(\mathcal{V},F,h)$ defined by   $\mathcal{V}= \mathbb{R}^N $ and
\begin{empheq}[left={\empheqlbrace}]{align}
F(\mathbf{x},z)&:=A \mathbf{x}+ \mathbf{C}z,\label{linear state eq}\\
h(\mathbf{x})&:= W \mathbf{x}. \label{linear readout eq}
\end{empheq}
\begin{description}
\item [(i)] The state equation associated to \eqref{linear state eq} has a unique solution in $\ell^{\infty}_-(\mathbb{R}^N)$ for each input in $\ell^{\infty}_-(\mathbb{R}) $ (we call this property the $\left(\ell^{\infty}_-(\mathbb{R}^N), \ell^{\infty}_-(\mathbb{R})\right)$-ESP) if and only if $\rho(A)< 1$, where $\rho(A) $ stands for the spectral radius of $A$.
\item [(ii)] In the remainder of this proposition suppose that $\rho(A)< 1$. Then, there exists a state filter $U ^F: \ell^{\infty}_-(\mathbb{R}) \longrightarrow \ell^{\infty}_-(\mathbb{R}^N) $ and a corresponding state-space filter $U ^F_h: \ell^{\infty}_-(\mathbb{R}) \longrightarrow \ell^{\infty}_-(\mathbb{R}) $ given by
\begin{equation}
\label{filter linear state-space}
U ^F({\bf z})_t:=\sum_{j=0}^{\infty}A ^j \mathbf{C}z_{t-j} \quad \mbox{and} \quad U ^F_h({\bf z})_t:=W\sum_{j=0}^{\infty}A ^j \mathbf{C}z_{t-j}, \quad \mbox{respectively.}
\end{equation}
The state-space filter $U ^F_h $ has the input forgetting property.
\item [(iii)] The set $\mathcal{V}_R\subset \mathcal{V} $ of reachable states defined in \eqref{set of reachable states} of $(\mathcal{V},F,h)$ is given by
\begin{equation}
\label{reachable for linear}
\mathcal{V}_R= {\rm span} \left\{\mathbf{C},A \mathbf{C}, A ^2 \mathbf{C}, \ldots, A^{N-1}\mathbf{C}\right\}.
\end{equation}
\item [(iv)] Given $\mathbf{x} \in \mathcal{V}_R  $, the set of indistinguishable states of $\mathbf{x} $ in $\mathcal{V}_R  $ is given by the coset
\begin{equation}
\label{indistinguishable set}
I_{F,h}^{\mathbf{x}}:= \mathbf{x}+I_{F,h} \quad \mbox{with} \quad I_{F,h}:=\bigcap_{i=0}^{N-1}\ker WA ^i.
\end{equation}
The state-space system $(\mathcal{V},F,h)$ is hence observable if and only if $I_{F,h}= \left\{{\bf 0}\right\} $. {This condition is equivalent to the to the maximality of the rank of the observability matrix $O(A,W)$ defined by
\begin{equation}
\label{observability statement}
O(A,W)=
\left(
\begin{array}{c}
W\\
WA\\
\vdots\\
WA^{N-1}
\end{array}
\right), 
\mbox{ that is, } \, {\rm rank}\, O(A,W)=N  \mbox{ if and only if } I_{F,h}= \left\{{\bf 0}\right\}.
\end{equation}}
\end{description}
\end{proposition}

\begin{remark}
\normalfont
The dimension of $\mathcal{V}_R $  in \eqref{reachable for linear} coincides with the rank of the {\bfi  controllability } or {\bfi  reachability matrix} $R(A, \mathbf{C})$ defined by
\begin{equation*}
R(A, \mathbf{C}):= \left(\mathbf{C}\mid A \mathbf{C}\mid \cdots\mid A^{N-1}\mathbf{C}\right).
\end{equation*}
When this rank is maximal, the linear system $(\mathcal{V},F,h) $ is strongly reachable in the sense of the Definition \ref{Reachable and Observable} and also in the control theoretical sense (see \cite{Kalman2010, sontag:book}). It has been shown in \cite{RC15} that if $A$ is diagonalizable then $R(A, \mathbf{C}) $ has maximal rank if and only if all the eigenvalues in the spectrum $\sigma(A) $ of $A$ are distinct and in the linear decomposition $\mathbf{C}=\sum_{i=1}^N c _i \mathbf{v} _i$, with $\left\{\mathbf{v}_1, \ldots, \mathbf{v}_N\right\}$ a basis of eigenvectors of $A$, all the coefficients $c _i $, with $i \in \left\{1, \ldots, N\right\} $, are non-zero.

The dimension of $\mathcal{V}_R $ also coincides with the so-called {\bfi  memory capacity} \cite{Jaeger:2002} of the recurrent network associated to \eqref{linear state eq}-\eqref{linear readout eq}. This fact has been recently proved in \cite{RC15}.
\end{remark}

\noindent\textbf{Proof of the Proposition.\ \ (i) and (ii)} We first show that if $\rho(A)\geq 1 $ then $(\mathcal{V},F) $ cannot have the $\left(\ell^{\infty}_-(\mathbb{R}^N), \ell^{\infty}_-(\mathbb{R})\right)$-ESP. Let $\lambda\geq 1 $ be one of the elements in the spectrum $\sigma (A)$ and let $\mathbf{v} \in \mathbb{R}^N $ be an associated norm-one eigenvector. Let $\mathbf{x} ^\lambda \in \ell^{\infty}_-(\mathbb{R}^N) $ be defined by $\mathbf{x} ^\lambda _t:= \lambda ^t \mathbf{v}  $, $t \in \mathbb{Z}_{-}  $. It is clear that as $\lambda\geq 1  $ then $\left\|\mathbf{x} ^\lambda\right\| _{\infty} = \left\|\mathbf{v}\right\|=1$. Moreover, $\mathbf{x}  ^\lambda $  is a solution of the system associated to $F$ with zero input because for any $t \in \mathbb{Z}_{-}  $ we have
\begin{equation*}
F(\mathbf{x}^\lambda _{t-1}, 0)=A \mathbf{x}^\lambda _{t-1}= \lambda^{t-1}A \mathbf{v}= \lambda ^t \mathbf{v}= \mathbf{x} ^\lambda _t.
\end{equation*}
Since ${\bf 0} \in \ell^{\infty}_-(\mathbb{R}^N) $ is also a solution for the same input, then $(\mathcal{V},F) $ does not have the ESP. What we just proved is equivalent to stating that if $(\mathcal{V},F) $ has the $\left(\ell^{\infty}_-(\mathbb{R}^N), \ell^{\infty}_-(\mathbb{R})\right)$-ESP then $\rho< 1  $ necessarily. 

Conversely, suppose that $\rho(A)<1 $. We now show that first, for any ${\bf z}\in \ell^{\infty}_-(\mathbb{R})  $ the sequence $\mathbf{x} \in \ell^{\infty}_-(\mathbb{R}^N) $ whose terms $\mathbf{x} _t $ are defined by 
\begin{equation}
\label{form of solution}
\mathbf{x} _t:=\sum_{j=0} ^{\infty} A ^j \mathbf{C} z_{t-j}
\end{equation}
is a solution of $(\mathcal{V}, F) $ for the input ${\bf z} $ and second, that this solution is unique. In order to show that \eqref{form of solution} is a solution,  we first recall that by Gelfand's formula (see \cite{lax:functional:analysis}) $\lim\limits_{k \rightarrow \infty}\vertiii{A ^k}^{1/k}=\rho(A)<1$, which implies the existence of a number $k _0 \in \mathbb{N}$ such that $\vertiii{A ^k}<1 $, for all $k\geq k _0 $. Consequently, the infinite sum
\begin{equation}
\label{sum converges factor}
\sum_{j=0}^{\infty}A ^j= \mathbb{I}_N+A+ \cdots+A^{k _0-1}+\sum _{j=1}^{\infty}\sum_{i=0}^{k _0-1}A^{j k _0}A ^i= \sum _{j=0}^{\infty}\sum_{i=0}^{k _0-1}A^{j k _0}A ^i
\end{equation}
converges in operator norm because as $\vertiii{A^{jk _0}A ^i}\leq \vertiii{A^{k _0}}^j\vertiii{A ^{i}} $ for all $j \in \mathbb{N}  $, $i \in  \left\{0, \ldots, k _0-1\right\}$ then \eqref{sum converges factor} implies that 
\begin{equation}
\label{rho smaller 1 implies finite l1norm}
\vertiii{\sum_{j=0}^{\infty}A ^j}\leq \sum_{i=0}^{k _0-1} \frac{\vertiii{A ^i}}{1-\vertiii{A^{k _0}}}< \infty.
\end{equation}
This inequality,  \eqref{sum converges factor},  and \eqref{form of solution} imply that 
\begin{equation*}
\left\|\mathbf{x} _t\right\|\leq \left\|\sum _{j=0}^{\infty}\sum_{i=0}^{k _0-1}A^{j k _0}A ^i \mathbf{C}z_{t-(jk _0+i)}  \right\|\leq 
\left(\sum _{j=0}^{\infty}\sum_{i=0}^{k _0-1}\vertiii{A^{k _0}}^j \vertiii{A ^i}  \right)\left\|\mathbf{C}\right\| \left\|{\bf z}\right\|_{\infty}=
 \sum_{i=0}^{k _0-1} \frac{\vertiii{A ^i}}{1-\vertiii{A^{k _0}}}\left\|\mathbf{C}\right\| \left\|{\bf z}\right\|_{\infty},
\end{equation*}
which shows that the series in \eqref{form of solution} are convergent and also that 
\begin{equation*}
\left\|\mathbf{x}\right\|_{\infty}\leq  \sum_{i=0}^{k _0-1} \frac{\vertiii{A ^i}}{1-\vertiii{A^{k _0}}}\left\|\mathbf{C}\right\| \left\|{\bf z}\right\|_{\infty}< \infty.
\end{equation*}
The fact that $\mathbf{x} \in \ell^{\infty}_-(\mathbb{R}^N) $ is a solution of $(\mathcal{V},F ) $ for the input ${\bf z} \in \ell^{\infty}_-(\mathbb{R}) $ is a straightforward verification. Suppose now that $\overline{\mathbf{x}} \in \ell^{\infty}_-(\mathbb{R}^N) $ is another solution of $(\mathcal{V},F ) $ for the same input, that is, $\overline{\mathbf{x}} _t=A \overline{\mathbf{x}} _{t-1}+ \mathbf{C} z _t$, for all $t \in \mathbb{Z}_{-} $. This implies that $\mathbf{x}- \overline{\mathbf{x}} \in \ell^{\infty}_-(\mathbb{R}^N) $ is a solution of $(\mathcal{V},F ) $ for the zero input and hence 
\begin{equation}
\label{to be iterated}
\mathbf{x} _t- \overline{\mathbf{x}} _t=A \left(\mathbf{x} _{t-1}- \overline{\mathbf{x}} _{t-1}\right), \quad \mbox{for all} \quad t \in \mathbb{Z}_{-}.
\end{equation}
 Using the same decomposition as in  \eqref{sum converges factor}, we have that for any $l \in \mathbb{N}   $ there exists $j \in \mathbb{N} $ and $i \in \left\{1, \ldots, k _0-1\right\} $ such that $A ^l= A^{j k _0}A ^i $. Hence, by iterating \eqref{to be iterated} we have that $\mathbf{x} _t- \overline{ \mathbf{x}} _t=A ^l \left(\mathbf{x} _{t-l}- \overline{ \mathbf{x}} _{t-l}\right) $ and therefore
\begin{equation*}
\left\|\mathbf{x} _t- \overline{\mathbf{x}} _t\right\|\leq \vertiii{A^{jk _0}A ^i} \left\|\mathbf{x}_{t-l}- \overline{\mathbf{x}} _{t-l}\right\|\leq \vertiii{A^{k _0}}^j\vertiii{A ^i}\left\|\mathbf{x}- \overline{\mathbf{x}}\right\|_ \infty.
\end{equation*}
Taking the limit $j \rightarrow \infty $ in this inequality, we obtain that $\left\|\mathbf{x} _t- \overline{\mathbf{x}} _t\right\|=0 $, for all $t \in \mathbb{Z}_{-}  $, which guarantees that $\mathbf{x} = \overline{\mathbf{x}} $, as required.

Finally, we show that when $\rho(A)<1 $ then the filter $U^F_h $ in \eqref{filter linear state-space} has the input forgetting property. Notice first that \eqref{filter linear state-space} amounts to a convolution representation for $U^F_h$, that is, $U^F_h({\bf z})= \boldsymbol{\Psi} \ast {\bf z}$, for any ${\bf z} \in \ell_{-}^{\infty}(\mathbb{R}) $, where $\Psi_{-j}= WA ^j \mathbf{C} $, $j \in \mathbb{N}  $. If we show that $\boldsymbol{\Psi} \in \ell_{-}^{1}(\mathbb{R}) $, then an argument similar to \eqref{implies fmp with l1} proves that $U^F _h $ has the input forgetting property. This is the case because by \eqref{rho smaller 1 implies finite l1norm}
\begin{equation*}
\left\|\boldsymbol{\Psi}\right\|_1= \left\|\sum _{j=0}^{\infty} WA^j \mathbf{C}\right\|\leq \vertiii{W} \vertiii{\sum _{j=0}^{\infty} A^j} \left\|\mathbf{C}\right\|\leq 
\vertiii{W}  \left\|\mathbf{C}\right\| \sum_{i=0}^{k _0-1} \frac{\vertiii{A ^i}}{1-\vertiii{A^{k _0}}}< \infty.
\end{equation*}

\medskip

\noindent {\bf (iii)} First of all, since by \eqref{filter linear state-space} the state functional $H ^F: \ell^{\infty}_-(\mathbb{R})\longrightarrow \mathbb{R}^N $ is linear and given by $H^F({\bf z}):=\sum_{j=0}^{\infty}A ^j \mathbf{C}z_{-j} $, we can immediately conclude that the reachable set $\mathcal{V}_R= H ^F(\ell^{\infty}_-(\mathbb{R})) \subseteq \mathbb{R} ^N$ is a vector subspace of $\mathbb{R} ^N$. We now establish \eqref{reachable for linear} by double inclusion. The inclusion 
$\mathcal{V}_R\supseteq {\rm span} \left\{\mathbf{C},A \mathbf{C}, A ^2 \mathbf{C}, \ldots, A^{N-1}\mathbf{C}\right\}$ is proved by applying $H ^F $ to inputs of the form 
\begin{equation*}
{\bf z}^j:= (\ldots, 0,\underbrace{1}_{\mbox{(-$j$)th entry}},0, \ldots,0) \in \ell^{\infty}_-(\mathbb{R}), \quad \mbox{with} \quad j \in \left\{0, \ldots, N-1\right\}. 
\end{equation*}
Conversely, let $\mathcal{V} _R^l $ be the reachable set associated to the truncated functional $H^F_l({\bf z}):=\sum_{j=0}^{l}A ^j \mathbf{C}z_{-j} $, $l \in \mathbb{N}  $. It is obvious that $\mathcal{V} _R^{N-1} \subset {\rm span} \left\{\mathbf{C},A \mathbf{C}, A ^2 \mathbf{C}, \ldots, A^{N-1}\mathbf{C}\right\}$. We now prove by induction that  $\mathcal{V} _R^{N+i} \subset {\rm span} \left\{\mathbf{C},A \mathbf{C}, A ^2 \mathbf{C}, \ldots, A^{N-1}\mathbf{C}\right\}$,  for all $i \in \mathbb{N}  $. First, by the Cayley-Hamilton Theorem \cite[Theorem 2.4.3.2]{horn:matrix:analysis} there exist constants $\left\{\alpha_0, \ldots, \alpha_{N-1}\right\} $ not all zero such that 
\begin{equation}
\label{ham cay for A}
A ^N= \alpha _0 \mathbb{I}_N+ \alpha _1A+ \cdots + \alpha _{N-1}A^{N-1}
\end{equation}
and hence $H^F_N({\bf z}):=\sum _{i=0}^{N-1}\alpha _iA ^{i} \mathbf{C}z_{-N}+\sum_{j=0}^{N-1}A ^j \mathbf{C}z_{-j} $, which shows that $\mathcal{V} _R^{N} \subset {\rm span} \left\{\mathbf{C},A \mathbf{C}, A ^2 \mathbf{C}, \ldots, A^{N-1}\mathbf{C}\right\}$. In order to prove the induction step, suppose that the inequality $\mathcal{V} _R^{N+(i-1)} \subset {\rm span} \left\{\mathbf{C},A \mathbf{C}, A ^2 \mathbf{C}, \ldots, A^{N-1}\mathbf{C}\right\}$ holds for a certain $i\in \mathbb{N}$. Again, using \eqref{ham cay for A}, we have that
\begin{multline*}
H^F_{N+i}({\bf z})=A ^NA^i \mathbf{C}z_{-(N+i)}+\sum_{j=0}^{N+(i-1)}A^j \mathbf{C}z_{-j}=\sum _{l=0}^{N-1}\alpha _lA ^{l+i} \mathbf{C}z_{-(N+i)}+\sum_{j=0}^{N+(i-1)}A^j \mathbf{C}z_{-j}\\=\sum_{j=0}^{i-1}A ^j \mathbf{C}z_{-j}+\sum_{j=i}^{N+i-1}A^j \mathbf{C}(z_{-j}+ \alpha_{j-i}z_{-(N+i)})=H^F_{N+(i-1)}(\overline{{\bf z}}),
\end{multline*}
with 
\begin{equation*}
\overline{{\bf z}}= \left(\ldots, z_{-(N+i)}, z_{-(N+i-1)}+ \alpha_{N-1}z_{-(N+i)}, \ldots, z_{-i-1}+ \alpha_{1}z_{-(N+i)},z_{-i}+ \alpha_{0}z_{-(N+i)}, z_{-(i-1)}, \ldots, z_{-1}, z _0\right),
\end{equation*}
which shows that $\mathcal{V}_R^{N+i} \subset \mathcal{V}_R^{N+(i-1)}\subset {\rm span} \left\{\mathbf{C},A \mathbf{C}, A ^2 \mathbf{C}, \ldots, A^{N-1}\mathbf{C}\right\} $ and hence proves the induction step. This inclusion also implies that
\begin{equation*}
\mathcal{V}_R= \bigcup_{l=0}^{\infty}\mathcal{V} _R^{l}\subset {\rm span} \left\{\mathbf{C},A \mathbf{C}, A ^2 \mathbf{C}, \ldots, A^{N-1}\mathbf{C}\right\},
\end{equation*}
as required.

\medskip

\noindent {\bf (iv)} Let $\mathbf{x} ^1, \mathbf{x} ^2 \in \mathcal{V}_R $ be two indistinguishable states of $(\mathcal{V}, F, h)$. By definition, this implies that there exist ${\bf z}^1, {\bf z} ^2 \in \ell^{\infty}_-(\mathbb{R}) $ such that $\mathbf{x} ^1=H ^F({\bf z}^1) $, $\mathbf{x} ^2=H ^F({\bf z}^2) $,  and that for any $\widetilde{{\bf z}} \in \mathbb{R}^T $  and any $T \in \mathbb{N}  $ we have that $H^F_h({\bf z}^1\widetilde{{\bf z}} )=H^F_h({\bf z}^2\widetilde{{\bf z}} )$. By \eqref{filter linear state-space}, this is equivalent to
\begin{equation*}
W \left(\sum_{j=0}^{\infty}A^{j+T}\mathbf{C}z_{-j}^1+\sum_{j=0}^{T-1}A^j \mathbf{C}\widetilde{z}_{T-j}\right)=W \left(\sum_{j=0}^{\infty}A^{j+T}\mathbf{C}z_{-j}^2+\sum_{j=0}^{T-1}A^j \mathbf{C}\widetilde{z}_{T-j}\right),
\end{equation*}
which amounts to
\begin{equation*}
W A^T\sum_{j=0}^{\infty}A^{j}\mathbf{C}z_{-j}^1=W A^T \sum_{j=0}^{\infty}A^{j}\mathbf{C}z_{-j}^2,
\end{equation*}
and is in turn equivalent to the relation $WA^T(\mathbf{x} ^1- \mathbf{x} ^2)=0 $, for all $T \in \mathbb{N} $ or, analogously, to $\mathbf{x} ^1- \mathbf{x} ^2 \in \bigcap_{j=0}^{\infty}\ker WA ^j $. In order to conclude the proof, it hence suffices to show that
\begin{equation*}
\bigcap_{j=0}^{\infty}\ker WA ^j=\bigcap_{j=0}^{N-1}\ker WA ^j.
\end{equation*}
The inclusion $\bigcap_{j=0}^{\infty}\ker WA ^j\subseteq \bigcap_{j=0}^{N-1}\ker WA ^j $ is obvious. Conversely, we show by induction that 
\begin{equation}
\label{inclusion by induction}
\bigcap_{j=0}^{N-1}\ker WA ^j\subseteq \bigcap_{j=0}^{N-1+i}\ker WA ^j  \quad \mbox{for all $i \in \mathbb{N} $.} 
\end{equation}
The initialization step is proved using the Cayley-Hamilton Theorem as formulated in \eqref{ham cay for A}. Indeed:
\begin{equation*}
\bigcap_{j=0}^{N}\ker WA ^j=\ker \left(\sum _{j=0}^{N-1}\alpha_j WA^j \right)\bigcap \left(\bigcap_{j=0}^{N-1}\ker WA ^j\right),
\end{equation*}
which obviously implies that $\bigcap_{j=0}^{N-1}\ker WA ^j\subseteq \bigcap_{j=0}^{N}\ker WA ^j $. In order to prove the induction step, suppose that \eqref{inclusion by induction} holds for a given $i \in \mathbb{N} $. Given that 
\begin{equation*}
\bigcap_{j=0}^{N+i}\ker WA ^j= \left(\bigcap_{j=0}^{N-1+i}\ker WA ^j\right)\bigcap \ker WA ^{N+i},
\end{equation*}
by the induction hypothesis we just need to show that $\bigcap_{j=0}^{N-1}\ker WA ^j\subseteq \ker WA ^{N+i} $.  This inclusion is easily established by using again the Cayley-Hamilton Theorem, which implies that $\ker WA ^{N+i} = \ker \sum_{j=0}^{N-1} \alpha _iWA^{i+j} $. The inclusion then follows from the induction hypothesis. 

{Finally, the statement \eqref{observability statement} follows in a straightforward manner from observing that:
\begin{equation*}
\ker O(A,W)=I_{F,h}.\quad \blacksquare
\end{equation*}}

\medskip

Before we use Theorem \ref{Canonical realization of input/output systems} in order to show that fading memory linear filters admit a linear canonical state-space realization, we motivate that result with an elementary example that hints how such construction may be obtained.

\begin{example}
\normalfont
{\bf Canonical realization of finite-memory linear filters.} Consider the finite-memory linear filter
\begin{equation*}
U({\bf z})_t=\sum_{j=0}^{-N+1} \Psi _j z _{t+j},\quad \mbox{with} \quad {\mathbf{\Psi}}\in \mathbb{R}^{N}, {\bf z} \in  \ell^{\infty}_-(\mathbb{R}), t \in \mathbb{Z}_{-},
\end{equation*}
and some $N \in \mathbb{N} $. Using the definition of the Nerode equivalence $\sim_I $ on the input space introduced in the proof of Theorem \ref{Canonical realization of input/output systems}, it is easy to see that ${\bf z}^1, {\bf z} ^2\in \ell^{\infty}_-(\mathbb{R}) $ are such that  ${\bf z}^1\sim _I {\bf z} ^2  $ if and only if $\left({z}^1_{-N+1}, \ldots, {z} ^1_0\right)=\left({z}^2_{-N+1}, \ldots, {z} ^2_0\right) $ and hence $\ell^{\infty}_-(\mathbb{R})/\sim _I  $ can be identified in this case with $\mathbb{R}^N  $ via the map
\begin{equation*}
\begin{array}{cccc}
\phi: &\ell^{\infty}_-(\mathbb{R})/\sim _I   & \longrightarrow & \mathbb{R} ^N\\
	&[ {\bf z} ] &\longmapsto &\left({z}_{-N+1}, \ldots, {z}_0\right).
\end{array}
\end{equation*}
\end{example}
With this identification, $\ell^{\infty}_-(\mathbb{R})/\sim _I $ inherits the vector space structure of $\ell^{\infty}_-(\mathbb{R}) $ and, moreover, the canonical state-space realization \eqref{realizing state system} introduced in Theorem \ref{Canonical realization of input/output systems} is given by
\begin{empheq}[left={\empheqlbrace}]{align}
F \left(\left({z}_{-N+1}, \ldots, {z}_0\right), \widetilde{z}\right)&=\left({z}_{-N+2}, \ldots, {z}_0, \widetilde{z}\right),\notag\\
h \left(\left({z}_{-N+1}, \ldots, {z}_0\right)\right)&= \sum_{j=0}^{-N+1} \Psi _j z _{j},\notag
\end{empheq}
or, in matrix form:
\begin{empheq}[left={\empheqlbrace}]{align}
F \left(\left({z}_{-N+1}, \ldots, {z}_0\right), \widetilde{z}\right)&=\underbrace{\left(
\begin{array}{cccc}
0 & 1 &  & \\
	&\ddots & \ddots & \mbox{\Huge 0}\\
\mbox{\Huge 0}& & \ddots & 1\\
 & & &0
\end{array}
\right)}_{{A}}
\left(
\begin{array}{c}
{z}_{-N+1}\\
{z}_{-N+2}\\
\vdots\\
z _0
\end{array}
\right)+
\underbrace{
\left(
\begin{array}{c}
0\\
0\\
\vdots\\
1
\end{array}
\right)}_{\mathbf{C}}\widetilde{z}
,\label{definition of a and c}\\
h \left(\left({z}_{-N+1}, \ldots, {z}_0\right)\right)&=\left(\Psi_{-N+1}, \ldots, \Psi _0\right)
\left(
\begin{array}{c}
{z}_{-N+1}\\
{z}_{-N+2}\\
\vdots\\
z _0
\end{array}
\right).\notag
\end{empheq}
By  Theorem \ref{Canonical realization of input/output systems}, this realization of $U$ is canonical. An observation that will be key in the next result is that the nilpotent matrix $A$ in \eqref{definition of a and c} is the projection onto the quotient space $\ell^{\infty}_-(\mathbb{R})/\sim_{I} $ of the time delay operator $T_{-1}: \ell^{\infty}_-(\mathbb{R}) \longrightarrow \ell^{\infty}_-(\mathbb{R}) $ and that the input vector $ \mathbf{C} $ is a matrix expression for the projected version of the inclusion 
\begin{equation}
\label{definition i0}
\begin{array}{cccc}
i_0:& \mathbb{R} &\hookrightarrow &\ell^{\infty}_-(\mathbb{R}) \\
	&z &\longmapsto & (\ldots,0,z).
\end{array}
\end{equation}

\begin{theorem}[Canonical realization of linear fading memory filters]
\label{Canonicalization of linear filters}
Let $U: \ell^{\infty}_-(\mathbb{R}) \longrightarrow \ell^{\infty}_-(\mathbb{R}) $ be a linear, causal, and time-invariant filter such that the associated functional $H _U: \ell^{\infty}_-(\mathbb{R}) \longrightarrow \mathbb{R}  $ has the fading memory property. Then: 
\begin{description}
\item [(i)] The quotient space $\mathcal{V}:=\ell^{\infty}_-(\mathbb{R})/\sim_{I} $ has a natural vector space structure inherited from $\ell^{\infty}_-(\mathbb{R}) $. The time delay operator $T _{-1} $ and the inclusion in \eqref{definition i0} can be naturally projected to two linear maps $\mathcal{A}:=[T _{-1}] \in L(\mathcal{V}, \mathcal{V})$, $\mathcal{C}:=[i _0] \in L( \mathbb{R}, \mathcal{V}) $, as well as the functional $H _U $ that we use to define ${\cal W}:=[H _U] \in L(\mathcal{V}, \mathbb{R}) $.
\item [(ii)] The state-space system $(\mathcal{V}, \mathcal{F}, \mathpzc{h}) $ with $\mathcal{F}(\mathbf{v},z ):= \mathcal{A}( \mathbf{v})+ \mathcal{C} (z)  $ and $\mathpzc{h}(\mathbf{v}):= {\cal W}(\mathbf{v}) $ is a canonical linear realization of $U$.
\item [(iii)] Consider the action of the group $ {\rm GL}(\mathcal{V})$ of all the linear automorphisms of $\mathcal{V} $ and its action $\Phi$ on the product $\left(L( \mathcal{V}, \mathcal{V})\times L(\mathbb{R} , \mathcal{V}) \times L(\mathcal{V}, \mathbb{R})\right) $ via the map
\begin{equation*}
\begin{array}{cccc}
\Phi: &{\rm GL}(\mathcal{V}) \times \left(L( \mathcal{V}, \mathcal{V})\times L(\mathbb{R} , \mathcal{V}) \times L(\mathcal{V}, \mathbb{R})\right) & \longrightarrow & \left(L( \mathcal{V}, \mathcal{V})\times L(\mathbb{R} , \mathcal{V}) \times L(\mathcal{V}, \mathbb{R})\right)\\
 	 &\left(B, \left(\mathcal{A}, \mathcal{C}, {\cal W}\right)\right)&\longmapsto & \left(B \mathcal{A} B ^{-1}, B \mathcal{C}, {\cal W} B ^{-1}\right).
\end{array}
\end{equation*}
All the canonical representations of $U$are given by the orbit of the triple $(\mathcal{A}, \mathcal{C}, {\cal W}) $ introduced in part {\bf (i)} and hence the space of canonical representations is isomorphic to the homogeneous manifold ${\rm GL}(\mathcal{V})/ {\rm GL}(\mathcal{V})_{(\mathcal{A}, \mathcal{C}, {\cal W})} $, with ${\rm GL}(\mathcal{V})_{(\mathcal{A}, \mathcal{C}, {\cal W})} $ the isotropy subgroup of the element $(\mathcal{A}, \mathcal{C}, {\cal W}) $.
\item [(iv)] If the canonical realization in {\bf (ii)} is finite-dimensional, then there exists $N \in \mathbb{N} $ such that $\mathcal{V}\simeq \mathbb{R}^N  $, where this isomorphism is implemented by a choice of basis in $\mathcal{V} $. There are also matrices $A \in \mathbb{M} _N$, $\mathbf{C} \in \mathbb{R}^N$, $W \in \mathbb{M}_{1,N}$ that express in that basis $\mathcal{A} $, $\mathcal{C} $, and ${\cal W}$, respectively. Let $(\mathbb{R}^N, F, h)$ be the system corresponding to $(\mathcal{V}, \mathcal{F}, \mathpzc{h}) $ in that basis. Then:
\begin{description}
\item [(a)] $\rho(A)<1 $.
\item [(b)] The set of reachable states of $(\mathbb{R}^N, F, h)$ coincides with $\mathbb{R}^N= {\rm span}\left\{ \mathbf{C}, A \mathbf{C}, A ^2\mathbf{C}, \ldots, A^{N-1} \mathbf{C}\right\}$.
\item [(c)] $I_{F,h}:=\bigcap _{i=1}^{N-1}\ker W A ^{i}= \left\{{\bf 0}\right\}$.
\item [(d)] $U ({\bf z})_t=\sum_{j=0}^{\infty}WA^j \mathbf{C}z_{t-j}$, with ${\bf z} \in \ell^{\infty}_-(\mathbb{R}) $, for all $t \in \mathbb{Z}_{-} $.
\item [(e)] Let $\boldsymbol{\Psi} \in \ell_{-}^{1}(\mathbb{R})$ be the unique element such that  $U ({\bf z})= \boldsymbol{\Psi}\ast {\bf z} $ for any ${\bf z} \in \ell^{\infty}_-(\mathbb{R}) $. Then, $\Psi_{-j}=WA ^j \mathbf{C} $, for any $j \in \mathbb{N} $.
\end{description} 
\end{description}
\end{theorem}

\noindent\textbf{Proof.\ \ (i) and (ii)} Since by  Proposition \ref{properties fmp filter linear} the fading memory property implies the input forgetting property,  any linear filter that satisfies the hypotheses in the statement satisfies too those in Theorem \ref{Canonical realization of input/output systems} and consequently has a unique (up to system isomorphism) canonical state-space realization. We shall now study the realization introduced in the proof of that theorem and shall also see that it has the linear form stated in part {\bf (ii)}. First of all, recall that by Proposition \ref{properties fmp filter linear} there exists a unique element $\boldsymbol{\Psi} \in \ell_{-}^{1}(\mathbb{R}) $ such that  $U ({\bf z})= \boldsymbol{\Psi}\ast {\bf z} $, for any ${\bf z} \in \ell^{\infty}_-(\mathbb{R}) $. Using this convolution representation and the properties of infinite series it is obvious to prove that  if ${\bf z} ^1\sim_I \overline{{\bf z}}^1 $ and ${\bf z} ^2\sim_I \overline{{\bf z}}^2 $, then for any $\lambda \in \mathbb{R}$ we have that  $\lambda{\bf z} ^1+ {\bf z} ^2\sim_I \lambda\overline{{\bf z}}^1 + \overline{{\bf z}}^2$. This implies that the sum and multiplication by scalars in $\ell_{-}^{\infty}(\mathbb{R}) $ drop to the quotient space $\mathcal{V}:=\ell^{\infty}_-(\mathbb{R})/\sim_{I} $, making it into a vector space. 

Also, using the convolution representation of $U$ it is easy to prove that both the time delay operator $T _{-1} $, the inclusion in \eqref{definition i0}, and the functional $H _U $ can be naturally projected to the linear maps $\mathcal{A}:=[T _{-1}] \in L(\mathcal{V}, \mathcal{V})$, $\mathcal{C}:=[i _0] \in L( \mathbb{R}, \mathcal{V}) $, and  ${\cal W}:=[H _U] \in L(\mathcal{V}, \mathbb{R}) $, that are uniquely determined by the equalities:
\begin{eqnarray*}
\mathcal{A}\circ \pi_{\sim_I}&=&\pi_{\sim_I} \circ T _{-1},\\
\mathcal{C} &= &\pi_{\sim_I} \circ i _0,\\
{\cal W} \circ \pi_{\sim_I}&= & H _U,
\end{eqnarray*}
where  $\pi_{\sim_I}: \ell_{-}^{\infty}(\mathbb{R}) \longrightarrow \mathcal{V}:=\ell^{\infty}_-(\mathbb{R})/\sim_{I}$ is the canonical projection.

These maps can be used to rewrite the canonical realization proposed in \eqref{realizing state system} as
\begin{equation*}
F([\mathbf{z}], \widetilde{{ z}}):= [{\bf z} \widetilde{{ z}}]=\pi_{\sim_I} (T _{-1}({\bf z})+ i _0(\widetilde{{ z}}))= \mathcal{A}([{\bf z}])+ \mathcal{C}(\widetilde{z}) \quad \mbox{and} \quad
h([{\bf z}]):= H _U ({\bf z})= {\cal W}([{\bf z}]), 
\end{equation*}
as required.

\medskip

\noindent {\bf (iii)} is a consequence of Corollary \ref{uniqueness corollary} and the equalities \eqref{isomorphic state map}-\eqref{isomorphic readout map}. Finally, {\bf (iv)} is a corollary of the characterization in Proposition \ref{Linear state-space realizations with semi-infinite inputs}. \quad $\blacksquare$

\medskip

In the previous theorem we showed that as a Corollary of the Canonical Realization Theorem \ref{Canonical realization of input/output systems}, any fading memory linear filter admits a canonical linear state-space realization. We now show that the Canonicalization by Reduction Theorem \ref{Canonicalization by reduction} implies that any linear state-space system that has the echo state property and the fading memory property can be reduced to a canonical system that is also linear and has the same linear filter associated. The proof is a straightforward consequence of Theorem \ref{Canonicalization by reduction} and of Proposition \ref{Linear state-space realizations with semi-infinite inputs}.

\begin{theorem}[Canonicalization by reduction of linear state-space systems]
\label{Canonicalization by reduction of linear state-space systems}
Let $(\mathbb{R}^N, F, h) $ be the linear system determined by the maps $F(\mathbf{x},z):= A \mathbf{x}+ \mathbf{C}z $ and $h (\mathbf{x}):= { W} \mathbf{x} $, with $A \in \mathbb{M}_N $ such that  $\rho(A)<1 $, $\mathbf{C} \in \mathbb{R}^N  $, ${ W} \in \mathbb{M}_{1,N} $, and with inputs ${\bf z} \in \ell_{-}^{\infty}(\mathbb{R}) $. Denote by $U ^F_h: \ell^{\infty}_-(\mathbb{R}) \longrightarrow \ell^{\infty}_-(\mathbb{R}) $ the associated linear input forgetting filter given by \eqref{filter linear state-space}. Let $\mathcal{V}_R \subset \mathbb{R}^N $ and $I_{F,h}\subset \mathbb{R}^N $ be the subspaces defined in \eqref{reachable for linear} and \eqref{indistinguishable set}, respectively. Then $U ^F_h$ has a canonical linear realization $(\mathcal{V}, \overline{F}, \overline{h}) $ on the quotient vector space $\mathcal{V}:= \mathcal{V} _R/I_{F,h} $ given by the maps:
\begin{eqnarray}
\overline{F}([\mathbf{x}],z) &:= &[A]([\mathbf{x}])+[\mathbf{C}]z, \label{reduced linear 1}\\
\overline{h}([\mathbf{x}]) &:=&[W]([\mathbf{x}]),\label{reduced linear 2}
\end{eqnarray}
where if $\pi: \mathcal{V}_R \longrightarrow \mathcal{V}_R/I_{F,h} $ and $i:\mathcal{V}_R\hookrightarrow \mathbb{R} ^N $ are the canonical projection and inclusion, respectively, the linear maps $[A] \in L(\mathcal{V},\mathcal{V}), [W]\in L(\mathcal{V}, \mathbb{R}) $, and $   [C]\in \mathcal{V} $ in \eqref{reduced linear 1}-\eqref{reduced linear 2} uniquely determined by the relations
\begin{equation*}
[A] \circ \pi= \pi\circ A \circ i, \quad [\mathbf{C}]= \pi(\mathbf{C}), \quad \mbox{and} \quad [W] \circ \pi= \pi \circ W \circ i.
\end{equation*}
\end{theorem}

\section{Implicit reduction using RKHS}
\label{Implicit reduction using RKHS}

An important drawback of the dimension reduction  techniques proposed in the previous sections is the need to compute and characterize various reachable sets and quotient spaces, which may be complicated and hence may reduce the practical value of the results that we propose. A situation where these problems may be circumvented is the case is when the readout $h: {\cal X} \longrightarrow {\cal Y} $ in the observation equation  \eqref{rc readout eq} is linear. This situation is practically relevant since various state-space systems that satisfy this condition have been shown to exhibit universal approximation properties. It is the case, for instance, of state-affine systems \cite{RC6} and the widely used echo state networks \cite{RC7, RC20}.

The way we proceed in that setup consists in associating to any state-space system $F:  {\cal X}  \times {\cal Z} \longrightarrow {\cal X}$  that satisfies the echo state property, a reproducing kernel Hilbert space (RKHS) $\mathbb{H}$ (see, for instance, Chapter 6 in \cite{Mohri:learning:2012} or \cite{schoelkopf:smola:book} for a general presentation of kernel methods) using the state functional $H ^F: {\cal Z}^{\Bbb Z_-}\longrightarrow {\cal X}  $ as a feature map. We shall then show that when the state space ${\cal X}$  is a finite dimensional  Hilbert space, then $\mathbb{H}  $ is isometrically isomorphic to the linear span  $\overline{ {\cal X} _R} $ given by
\begin{equation}
\label{linear span def}
\overline{ {\cal X} _R} :={\rm span} \left\{{\cal X} _R\right\}= {\rm span}\left\{H ^F({\bf z})\mid {\bf z}\in  {\cal Z}^{\Bbb Z _-}\right\} 
\end{equation}
of the set of reachable states ${\cal X} _R $.

The importance of this characterization is in the fact that it allows us to show, using the classical Representer Theorem \cite[page 117]{Mohri:learning:2012}, that the search for an optimal readout with respect to the regularized empirical risk minimization associated to any loss can be reduced to the search for a readout defined on the smaller space $\overline{ {\cal X} _R} $ {\it without having to actually compute it}. We call this procedure {\bfi  implicit reduction}.

\paragraph{The RKHS associated to a state system.} Let $F:  {\cal X}  \times {\cal Z} \longrightarrow {\cal X}$ be a state equation such that the pair $\left({\cal X}, \langle \cdot , \cdot \rangle_{{\cal X}}\right) $ is a finite dimensional  Hilbert space and $F$ has the echo state property. Let $H ^F: {\cal Z}^{\Bbb Z_-}\longrightarrow {\cal X}  $ be the corresponding state functional. Define the {\bfi  kernel map} 
\begin{equation}
\label{definition kernel map}
\begin{array}{cccc}
K:&{\cal Z}^{\Bbb Z_-} \times {\cal Z}^{\Bbb Z_-} &\longrightarrow & \mathbb{R}\\
	&({\bf z}, {\bf z}')&\longmapsto & \langle H ^F({\bf z}) , H ^F({\bf z}') \rangle_{{\cal X}}.
\end{array}
\end{equation}
The map $K$ is obviously symmetric and positive semidefinite in the sense that for any $a _i \in \mathbb{R} $, ${\bf z} _i \in {\cal Z}^{\Bbb Z_-} $, $i \in \left\{1, \ldots, n\right\} $, we have that $\sum_{i,j=1}^n a _i a _jK({\bf z} _i, {\bf z} _j)\geq 0 $. Let $\left(\mathbb{H}, \langle \cdot , \cdot \rangle_{\mathbb{H}}\right)$ be the corresponding RKHS given by
\begin{equation}
\label{definition H for F}
\mathbb{H}:= \overline{{\rm span}\left\{K _{\bf z}:=K({\bf z}, \cdot ):{\cal Z}^{\Bbb Z_-} \longrightarrow \mathbb{R}\mid {\bf z} \in {\cal Z}^{\Bbb Z_-} \right\}}
\end{equation}
made out of finite linear combinations of elements of the type $K _{\bf z} $, ${\bf z} \in {\cal Z}^{\Bbb Z_-} $, together with all the limits of Cauchy sequences with respect to the metric induced by the inner product obtained as the linear extension of 
\begin{equation}
\label{definition hproduct}
\langle K _{\bf z}, K _{{\bf z}'}\rangle_{\mathbb{H}}=K({\bf z}, {\bf z}')=\langle H ^F({\bf z}) , H ^F({\bf z}') \rangle_{{\cal X}}, \quad {\bf z}, {\bf z}' \in {\cal Z}^{{\Bbb Z}_-}. 
\end{equation}
Note that in this setup, the reservoir functional $H ^F $ with respect to the kernel $K$ and the elements in $\mathbb{H} $ can be written as $K _{\bf z}(\cdot )= \langle H ^F({\bf z}) , H ^F(\cdot ) \rangle_{{\cal X}}$.

\begin{proposition}
Let $\left({\cal X}, \langle \cdot , \cdot \rangle_{{\cal X}}\right) $ be a finite dimensional  Hilbert space and let $F:  {\cal X}  \times {\cal Z} \longrightarrow {\cal X}$ be a state equation that satisfies the echo state property. Let  $\left(\mathbb{H}, \langle \cdot , \cdot \rangle_{\mathbb{H}}\right)$ be the associated RKHS introduced in \eqref{definition H for F}. Then
\begin{equation}
\label{first identity H}
\mathbb{H}= \left\{\langle {\bf W}, H ^F(\cdot )\rangle_{{\cal X}}\mid {\bf W} \in \overline{ {\cal X} _R}\right\}.
\end{equation}
Moreover, for any ${\bf W} _1, {\bf W} _2 \in \overline{ {\cal X} _R}$, we have that
\begin{equation}
\label{first identity H 2}
\left\langle \langle {\bf W}_1, H ^F(\cdot )\rangle_{{\cal X}},  \langle {\bf W}_2, H ^F(\cdot )\rangle_{{\cal X}}\right\rangle_{\mathbb{H}}= \langle {\bf W} _1, {\bf W} _2\rangle_{{\cal X}},
\end{equation}
and the map
\begin{equation}
\label{isomorphism H}
\begin{array}{cccc}
\Psi: &\left(\overline{{\cal X}_R}, \langle \cdot , \cdot \rangle_{{\cal X}}\right)& \longrightarrow & \left(\mathbb{H}, \langle \cdot , \cdot \rangle_{\mathbb{H}}\right)\\
	&{\bf W} &\longmapsto &\langle {\bf W}, H ^F(\cdot )\rangle_{{\cal X}}=:H_{\bf W}^F(\cdot )
\end{array}
\end{equation}
is an isometric isomorphism.
\end{proposition}

\noindent\textbf{Proof.\ \ } We first establish the identity \eqref{first identity H} by double inclusion. In order to show that ${\mathcal H}\subset \left\{\langle {\bf W}, H ^F(\cdot )\rangle_{{\cal X}}\mid W \in \overline{ {\cal X} _R}\right\} $ consider the element
\begin{equation}
\label{element type}
f=\sum _{i=1}^n a _iK_{{\bf z} _i} \in \mathbb{H},\quad \mbox{for some $a _1, \ldots, a _n\in \mathbb{R} $}. 
\end{equation}
Then $f (\cdot )=\sum_{i=1}^n \langle a _iH^F({\bf z}_i), H^F(\cdot ) \rangle_{{\cal X}}$. Hence, it is clear that if we set ${\bf W}:=\sum_{i=1}^n  a _iH^F({\bf z}_i)$ we can then obviously write that $f(\cdot )= \langle {\bf W}, H ^F(\cdot )\rangle_{{\cal X}} $, as required. More generally, what we just showed also proves that for any sequence $\left\{f _n\right\} _{n \in \mathbb{N}}$ of elements like \eqref{element type} there are elements ${\bf W}_n\in \overline{{\cal X}_R}  $ such that $f _n(\cdot )= \langle{\bf W}_n, H^F(\cdot )\rangle$. If we assume that $\left\{f _n\right\} _{n \in \mathbb{N}}$ is Cauchy then $\left\|f _n- f _m\right\|_{\mathbb{H}}\rightarrow 0 $ as $n,m \rightarrow\infty$. This in turn implies that for any ${\bf z} \in {\cal Z}^{\Bbb Z_-} $ we have that
\begin{equation*}
\left| f _n({\bf z})- f _m({\bf z}) \right|= \left| \langle K_{{\bf z}} , f _n- f _m \rangle_{\mathbb{H}} \right|\leq \left\|f _n- f _m\right\|_{\mathbb{H}}\left\|K_{{\bf z}}\right\|_{\mathbb{H}}\longrightarrow 0 \quad \mbox{as $n,m \rightarrow\infty$,}
\end{equation*}
which guarantees that $\left| \langle{\bf W}_n- {\bf W}_m, H^F({\bf z})\rangle \right| \rightarrow 0 $  and hence that $\left| \langle{\bf W}_n- {\bf W}_m, \mathbf{v}\rangle \right| \rightarrow 0 $ as $n,m \rightarrow \infty $, for any $\mathbf{v} \in \overline{{\cal X}_R} $. Now, since any vector $\mathbf{x} \in {\cal X} $ can be uniquely decomposed as $\mathbf{x}= \mathbf{v}+ \mathbf{v} ^\perp $ with  $\mathbf{v} \in \overline{{\cal X}_R} $ and $\mathbf{v}^\perp \in \overline{{\cal X}_R} ^\perp$, we also have that
\begin{equation*}
\left| \langle{\bf W}_n- {\bf W}_m, \mathbf{x}\rangle \right|= \left| \langle{\bf W}_n- {\bf W}_m, \mathbf{v}\rangle \right| \rightarrow 0 \quad \mbox{as $n,m \rightarrow \infty $, for any $\mathbf{x} \in \overline{{\cal X}_R} $.}
\end{equation*}
Given that ${\cal X} $ is finite dimensional, we can conclude that weak and strong convergence coincide and hence that $\left\|{\bf W}_n- {\bf W}_m\right\|_{{\cal X}} \rightarrow 0 $ as $n,m \rightarrow \infty $. Since ${\cal X} $ is complete then so is $\overline{{\cal X}_R} $ and hence there exists ${\bf W} \in \overline{{\cal X}_R} $ such that $\lim\limits _{n \rightarrow \infty} {\bf W}_n= {\bf W} $. It is easy to see that this implies that 
\begin{equation*}
f (\cdot ):=\lim\limits_{n \rightarrow \infty} f _n(\cdot )= \langle{\bf W}, H ^F (\cdot )\rangle_{{\cal X}},
\end{equation*}
as required. In order to prove the converse inclusion, note first that by definition, for any ${\bf W} \in \overline{ {\cal X}_R} $ there exist ${\bf z} _1, \ldots, {\bf z} _n \in {\cal Z}^{{\Bbb Z}_-} $ and $a _1, \ldots, a _n \in \mathbb{R}  $ such that ${\bf W}=\sum_{i=1}^n  a _iH^F({\bf z}_i)$. It is hence easy to see that
\begin{equation*}
\langle{\bf W}, H^F(\cdot )\rangle_{{\cal X}}= \sum_{i=1}^n  a _i\langle H^F({\bf z}_i), H^F(\cdot )\rangle_{{\cal X}}=\sum_{i=1}^n  a _iK_{{\bf z} _i}(\cdot ),
\end{equation*}
which is an element in $\mathbb{H}$, as required.

We now show the identity \eqref{first identity H 2}. Let ${\bf W} _1, {\bf W} _2 \in \overline{ {\cal X} _R}$ and let ${\bf W} _1=\sum_{i=1}^n  a _i^1H^F({\bf z}^1_i)$ and ${\bf W} _2=\sum_{i=1}^n  a _i^2H^F({\bf z}^1_2)$ two representations of the two vectors according to the definition of $\overline{ {\cal X}_R} $. Then, it is easy to see that
\begin{multline*}
\left\langle \langle {\bf W}_1, H ^F(\cdot )\rangle_{{\cal X}},  \langle {\bf W}_2, H ^F(\cdot )\rangle_{{\cal X}}\right\rangle_{\mathbb{H}}=\sum_{i=1}^{n _1}\sum_{j=1}^{n _2}a _i^1a _j^2\left\langle \langle H^F({\bf z}^1_i), H ^F(\cdot )\rangle_{{\cal X}},  \langle H^F({\bf z}^2_i), H ^F(\cdot )\rangle_{{\cal X}}\right\rangle_{\mathbb{H}}\\
=\sum_{i=1}^{n _1}\sum_{j=1}^{n _2}a _i^1a _j^2\left\langle K_{{\bf z} _i^1}, K_{{\bf z} _j^2}\right\rangle_{\mathbb{H}}=\sum_{i=1}^{n _1}\sum_{j=1}^{n _2}a _i^1a _j^2\left\langle H^F({\bf z}^1_i), H^F({\bf z}^2_j)\right\rangle_{\mathbb{X}}
= \langle {\bf W} _1, {\bf W} _2\rangle_{{\cal X}}.
\end{multline*}

Finally, we show that the map $\Psi $ in \eqref{isomorphism H} is an isometric isomorphism. First, it is clear that the map is linear, the equality \eqref{first identity H} guarantees that $\Psi $ is onto, and \eqref{first identity H 2} that it is an isometry. In order to show injectivity, suppose that ${\bf W}  \in \overline{ {\cal X} _R}$ is such that $\Psi({\bf W})(\cdot )=\langle {\bf W}, H ^F(\cdot )\rangle_{{\cal X}}=0 $. If we use a representation for ${\bf W}  $ of the type ${\bf W} =\sum_{i=1}^n  a _iH^F({\bf z}_i)$ we can write that
\begin{equation*}
\langle{\bf W}, {\bf W}\rangle_{{\cal X}}=\left\langle{\bf W}, \sum_{i=1}^n  a _iH^F({\bf z}_i)\right\rangle_{{\cal X}}=  \sum_{i=1}^n  a _i\Psi({\bf W})({\bf z}_i)=0,
\end{equation*}
which guarantees that ${\bf W}=0 $, as required. \quad $\blacksquare$

\paragraph{Estimation of the empirical risk minimizing readout.} A common estimation problem that appears when using in practice systems of the form \eqref{rc state eq}-\eqref{rc readout eq} and where the readout is linear is finding the readout vector ${\bf W} \in {\cal X} $ that minimizes the empirical risk associated to a prescribed loss function $L: {\cal Y}\times {\cal Y} \longrightarrow\mathbb{R}$ with respect to a finite sample of input/output observations. This is typically how one proceeds in reservoir computing (see the introduction section) where the state equation is fixed and only a linear observation equation is subjected to training. In that particular case and if a quadratic loss is used, the estimation problem reduces itself to a (eventually regularized) regression problem with as many covariates as the dimension of the state space ${\cal X} $, which is in most cases very large. It is in this context that for quadratic or more general losses, the possibility of reducing the dimensionality of the estimation problem to the dimension of $\overline{{\cal X}_R} $ using the RKHS technology that we just introduced may prove computationally advantageous. 

To be more specific, in the next proposition we will show two main fact as a consequence of the RKHS formulation of the estimation problem. First, that even though the optimization problem that provides the optimal readout is originally formulated in the space $\mathcal{X} $, {\it it can be reduced to the dimensionally smaller} $\overline{{\cal X}_R} $. Second, the Representer Theorem \cite[page 117]{Mohri:learning:2012} shows that {\it the optimal readout is in the ``span of the data"}; this is the well-known ``kernelization trick" that in our case is computational relevant in the presence of state spaces of dimension larger than the sample size. An important observation is that this second result {\it yields automatically a solution in the span $\overline{{\cal X}_R} $ of the reachable set ${\cal X}_R $ without actually having to compute it}.

We now introduce the different elements that are necessary for the statement of the Proposition. First, we will assume that the output space ${\cal Y} $  is one-dimensional, the state system $F: {\cal X} \times {\cal Y} \longrightarrow\mathcal{X} $ is fixed and satisfies the ESP, and we are provided with a finite sample 
$\{( {\bf Z}_{-i}, Y_{-i})\}_{i\in \left\{ 0, \dots, n-1\right\}}$ of size $n$ of input/output observations. For each time step $i\in \{ 0, \dots, n-1\}$ we define the truncated training sample for the input stochastic process $ {\bf Z}$ as
\begin{equation*} 
{\bf Z}_{-i}^{-n+1} := (\ldots,{\bf 0},{\bf 0},{\bf Z}_{-n+1},\ldots,{\bf Z}_{-i-1},{\bf Z}_{-i}),
\end{equation*}
that we use to define the {\bfi  training error} or the {\bfi  empirical risk} $\widehat{R}_n\left(H^F_{{\bf W}}\right)$ associated to the loss $L: {\cal Y}\times {\cal Y} \longrightarrow\mathbb{R}$ for the system $H^F_{{\bf W}} (\cdot )= \langle{\bf W}, H^F(\cdot )\rangle_{{\cal X}}$ with readout vector ${\bf W} \in {\cal X} $ as
\begin{equation*}
\widehat{R}_n\left(H^F_{{\bf W}}\right) = \frac{1}{n} \sum_{i=0}^{n-1} L\left( \left\langle{\bf W}, H^F({\bf Z}_{-i}^{-n+1})\right\rangle_{{\cal X}},{\bf Y}_{-i}\right).
\end{equation*}

\begin{proposition}
Let $F: {\cal X} \times {\cal Y} \longrightarrow\mathcal{X} $ be a state system that satisfies the ESP and let $L: {\cal Y} \times {\cal Y}\longrightarrow \mathbb{R}$ be a loss function with respect to the one-dimensional output space ${\cal Y} $. Let $\{( {\bf Z}_{-i}, Y_{-i})\}_{i\in \left\{ 0, \dots, n-1\right\}}$ be a sample of size $n$ of input/output observations. Let $\Omega: \mathbb{R}^+ \longrightarrow \mathbb{R} $ be a strictly increasing function. Then:
\begin{eqnarray}
\min_{W \in {\cal X}} \left\{\widehat{R}_n\left(H^F_{{\bf W}}\right)+ \Omega\left(\left\|{\bf W}\right\|_{{\cal X}}^2\right)\right\}&=&
\min_{W \in \overline{{\cal X}_R}} \left\{\widehat{R}_n\left(H^F_{{\bf W}}\right)+ \Omega\left(\left\|{\bf W}\right\|_{{\cal X}}^2\right)\right\}\label{first terms}\\
	&=&
\min_{H^F_{{\bf W}} \in \mathbb{H}} \left\{\widehat{R}_n\left(H^F_{{\bf W}}\right)+ \Omega\left(\left\|H_F^{{\bf W}}\right\|_{\mathbb{H} }^2\right)\right\},\label{last term prop}
\end{eqnarray} 
where $\mathbb{H}  $ is the RKHS introduced in \eqref{definition H for F}. The minimum in \eqref{last term prop} is realized by an element in $\mathbb{H} $ of the form 
\begin{equation}
\label{minimizer prop}
\sum_{i=0}^{n-1} \alpha _iK_{{\bf Z}_{-i}^{-n+1}}= \left\langle \sum_{i=0}^{n-1} \alpha _i H ^F \left({\bf Z}_{-i}^{-n+1}\right), H^F \left(\cdot \right)
\right\rangle_{{\cal X}}, \quad \mbox{for some} \quad \alpha_0, \alpha _1, \ldots, \alpha_{n-1} \in \mathbb{R}.
\end{equation}
The element $\sum_{i=0}^{n-1} \alpha _i H ^F \left({\bf Z}_{-i}^{-n+1}\right) \in \overline{{\cal X} _R} $ is the minimizer of the terms in \eqref{first terms}.
\end{proposition}

\noindent\textbf{Proof.\ \ } Given that any ${\bf W} \in {\cal X} $ can be uniquely decomposed as ${\bf W} = {\bf W}_{R}+{\bf W}_{R}^\perp $ with ${\bf W}_{R} \in {\cal X}_{R} $  and ${\bf W}_{R}^\perp \in {\cal X}_{R}^\perp $, we can write that
\begin{multline*}
\widehat{R}_n\left(H^F_{{\bf W}}\right)+ \Omega\left(\left\|{\bf W}\right\|_{{\cal X}}^2\right)=\widehat{R}_n\left(H^F_{{\bf W}_R}\right)+ \Omega\left(\left\|{\bf W}_R\right\|_{{\cal X}}^2+ \left\|{\bf W}_R^\perp\right\|_{{\cal X}}^2\right)\geq 
\widehat{R}_n\left(H^F_{{\bf W}_R}\right)+ \Omega\left(\left\|{\bf W}_R\right\|_{{\cal X}}^2\right),
\end{multline*}
where in the last inequality we used that, by hypothesis, $\Omega  $ is strictly increasing. This inequality implies that
\begin{equation*}
\min_{W \in {\cal X}} \left\{\widehat{R}_n\left(H^F_{{\bf W}}\right)+ \Omega\left(\left\|{\bf W}\right\|_{{\cal X}}^2\right)\right\}\geq 
\min_{W \in \overline{{\cal X}_R}} \left\{\widehat{R}_n\left(H^F_{{\bf W}}\right)+ \Omega\left(\left\|{\bf W}\right\|_{{\cal X}}^2\right)\right\}.
\end{equation*}
However, given that $\overline{{\cal X}_R} \subseteq {\cal X} $ the converse inequality also obviously holds, which proves the equality \eqref{first terms}.  The relation \eqref{last term prop} is a consequence of \eqref{first identity H} and also of the fact that by \eqref{first identity H 2}
\begin{equation*}
\left\|H_F^{{\bf W}}\right\|_{\mathbb{H} }^2=\left\langle \langle {\bf W}, H ^F(\cdot )\rangle_{{\cal X}},  \langle {\bf W}, H ^F(\cdot )\rangle_{{\cal X}}\right\rangle_{\mathbb{H}}= \left\|{\bf W}\right\|_{{\cal X}}^2.
\end{equation*}
Finally, the statement \eqref{minimizer prop} is a straightforward consequence of the Representer Theorem \cite[page 117]{Mohri:learning:2012}. \quad $\blacksquare$
\color{black}

\section{Conclusions}
\label{Conclusions}

In this paper we have extended the classical notion of canonical state-space realization to accommodate semi-infinite inputs so that it can be used as a dimension reduction tool in the framework of recurrent networks. We have formulated two main results that identify the so-called input forgetting property (introduced in Definition \ref{definition input forgetting property}) as the key hypothesis that guarantees the existence and uniqueness (up to system isomorphisms) of canonical realizations for causal and time-invariant input/output systems with semi-infinite inputs. 

The first result (Theorem \ref{Canonical realization of input/output systems}) shows that any causal and time-invariant filter with semi-infinite inputs that has the input forgetting property admits a canonical state-space realization that is unique up to system isomorphisms. The second one (Theorem \ref{Canonicalization by reduction}) uses a reduction approach similar to the one introduced in \cite{optimal:mm, reduction:optimal:cras} in the context of symmetric Hamiltonian systems to construct a canonical realization for a state-space system that has the input forgetting property system by using an ``optimally reduced" version of it, in the sense of those references.
These two results have been illustrated and applied in detail in Section \ref{Realization and canonicalization of linear filters} in the context of linear fading memory filters. 

The contributions in this paper should be considered just as a first step in the full understanding of this problem as, in comparison with the classical theory of forward-looking input-driven state-space systems, there are many deficiencies in the level of comprehension of several important mathematical issues. We now list a few of them that are part of our research agenda and that will be studied in forthcoming works:
\begin{itemize}
\item {\bfi  The geometric nature of  reachable sets by semi-infinite inputs} (see the definition in \eqref{set of reachable states}). Reachable sets are central objects in the context of continuous-time forward looking systems in connection with the notion of controllability (see \cite{sontag:book, Lewis2002, Bullo2005, bloch2015nonholonomicbook} and references therein). From the geometric viewpoint, this important application question has given rise to the notions of generalized foliation and distribution \cite{stefan, stefanb, sussmann, kms93}. Some of these results have a discrete-time counterpart (see, for instance, \cite{fliess1981group, jakubczyk1990controllability, kalman:Festschrift}) but the situation is mostly unknown when it comes to semi-infinite inputs. Some partial information \cite{manjunath:tino:jaeger} can be obtained by using the recent theory of nonautonomous dynamical systems \cite{Kloeden:Rasmussen}.
\item {\bfi  The geometric nature of the canonical state spaces obtained by reduction} (see the definition in Theorem \ref{Canonicalization by reduction}). Again, in other contexts like the reduction of symmetric Hamiltonian systems or control systems, this is a very well studied question (see \cite{Ortega2004, Marsden2007} for the autonomous case or \cite{van1981symmetries, nijmeijer1982controlled,  bloch2015nonholonomicbook} for the control case). The semi-infinite inputs framework presents new mathematical challenges that need to be addressed with innovative tools.
\item {\bfi  The geometric nature of the canonical realization state-spaces in Theorem \ref{Canonical realization of input/output systems}}. In the linear case treated in Section \ref{Realization and canonicalization of linear filters} we were able to easily  pinpoint the vector  space structure of the quotient space ${\cal X}:= {\cal Z}^{\mathbb{Z}_{-}}/\sim_{I} $ and to comfortably work with it. In more general nonlinear situations it is very difficult to answer even elementary questions (like the dimension) about the canonical state-space ${\cal X} $ even when we impose strong regularity assumptions on the original input space ${\cal Z}^{\mathbb{Z}_{-}} $.
\item {Even in the linear case, there is, as far as we know, no {\bf readily usable characterization of the situations in which the canonical realizations introduced in Theorem \ref{Canonicalization of linear filters} are finite dimensional}. Such criterion is necessary for the practical implementation of this result.}
\end{itemize}

\bigskip

\noindent {\bf Acknowledgments:} JPO acknowledges partial financial support  coming from the Research Commission of the Universit\"at Sankt Gallen and the Swiss National Science Foundation (grant number 200021\_175801/1). The authors thank the hospitality and the generosity of the FIM at ETH Zurich and the Division of Mathematical Sciences of the Nanyang Technological University, Singapore, where a significant portion of the results in this paper was obtained.

\noindent
\addcontentsline{toc}{section}{Bibliography}
\bibliographystyle{wmaainf}
%\bibliography{/Users/JPO/Dropbox/Public/GOLibrary}
%\bibliography{/Users/Lyudmila/Dropbox/Mendeley/GOLibrary}
%\bibliography{/Volumes/Staff/Mila/Library/BibTex/GOLibrary}

\end{document}